\documentclass[10pt,twoside,reqno]{tran-l}
\usepackage{graphicx}
\usepackage{amsmath}  	% \text{} etc
\usepackage{latexsym}	% LaTeX fonts
\usepackage{amssymb}  	% for \feedback symbol; superset of amsfonts
\usepackage{amsaddr}      % used to get proper journal addresses
\usepackage{amsthm}   	% AMS style theorem and proof environments
\usepackage{float}            % Allows figures to not float to next page
\usepackage{mathrsfs}	% special script style \mathscr{S} for \Set
\usepackage[all]{xy}	% Xy-Pic
\usepackage{mathtools}

\usepackage[T1,hyphens]{url}
\usepackage{pdflscape}
\usepackage{multirow}
\usepackage{pifont}          %%Used to do ZaphDingbats symbols
\usepackage{changebar}
\usepackage{fancyhdr}
\usepackage{ecltree}
\usepackage{placeins}
\usepackage{relsize,exscale}
\usepackage{bm}
\usepackage{chngcntr}
\usepackage[x11names]{xcolor}
\usepackage[colorlinks=true,linkcolor=blue,citecolor=SpringGreen4,urlcolor=blue]{hyperref}
\usepackage{rotating}
\usepackage{comment}
\usepackage{enumerate}

\usepackage{tikz-cd}

\usepackage[noend]{algpseudocode}
\usepackage{etoolbox}
\usepackage[utf8]{inputenc}
\usepackage[ruled, lined, longend, linesnumbered]{algorithm2e}

\makeatletter
\renewcommand{\tocsection}[3]{%
  \indentlabel{\@ifnotempty{#2}{\bfseries\ignorespaces#1 #2\quad}}\bfseries#3}
\renewcommand{\tocsubsection}[3]{%
  \indentlabel{\@ifnotempty{#2}{\ignorespaces#1 #2\quad}}#3}
\newcommand\@dotsep{4.5}
\def\@tocline#1#2#3#4#5#6#7{\relax
  \ifnum #1>\c@tocdepth % then omit
  \else
    \par \addpenalty\@secpenalty\addvspace{#2}%
    \begingroup \hyphenpenalty\@M
    \@ifempty{#4}{%
      \@tempdima\csname r@tocindent\number#1\endcsname\relax
    }{%
      \@tempdima#4\relax
    }%
    \parindent\z@ \leftskip#3\relax \advance\leftskip\@tempdima\relax
    \rightskip\@pnumwidth plus1em \parfillskip-\@pnumwidth
    #5\leavevmode\hskip-\@tempdima{#6}\nobreak
    \leaders\hbox{$\m@th\mkern \@dotsep mu\hbox{.}\mkern \@dotsep mu$}\hfill
    \nobreak
    \hbox to\@pnumwidth{\@tocpagenum{\ifnum#1=1\bfseries\fi#7}}\par% <-- \bfseries for \section page
    \nobreak
    \endgroup
  \fi}
\AtBeginDocument{%
\expandafter\renewcommand\csname r@tocindent0\endcsname{0pt}
}
\def\l@subsection{\@tocline{2}{0pt}{2.5pc}{5pc}{}}
\makeatother
\newtheorem{theorem}{Theorem}[section]

\theoremstyle{definition}
\newtheorem{definition}[theorem]{Definition}

\newtheorem{remark}[theorem]{Remark}
\numberwithin{equation}{section}

\usepackage{xcolor}
\usepackage{hyperref}      % Allows hyperlinks
\usepackage{colortbl}

\usepackage{booktabs}
\usepackage{xcolor}
\usepackage{array}
\usepackage{colortbl}

\renewcommand{\leq}{\leqslant}

\renewcommand{\geq}{\geqslant}

\usepackage{pgfplots}
\pgfplotsset{width=10cm,compat=1.9}
\usepackage{pgfplots,pgfplotstable}
\usepgfplotslibrary{colormaps}
\pgfplotsset{compat=newest}

\DeclareMathOperator{\sgn}{sgn}

\usepackage{xstring}
\def\typesystem#1{%
    \begingroup\expandarg
    \baselineskip=1.5\baselineskip% 1.5 to enlarge vertical space between lines
    \StrSubstitute{\noexpand#1}+{&+&}[\tempsystem]%
    \StrSubstitute\tempsystem={&=&}[\tempsystem]%
    \StrSubstitute\tempsystem,{\noexpand\cr}[\tempsystem]%
    $\vcenter{\halign{&$\hfil\strut##$&${}##{}$\cr\tempsystem\crcr}}$%
    \endgroup
}

\makeatletter
\newcommand*\bigcdot{\mathpalette\bigcdot@{.5}}
\newcommand*\bigcdot@[2]{\mathbin{\vcenter{\hbox{\scalebox{#2}{$\m@th#1\bullet$}}}}}
\makeatother

\usepackage[square, numbers]{natbib}  %<----------Make sure that natbib should go before apalike 9better still, take it out and reference it at the very top
\begin{document}

\title[Hidden Geometry of BGC Stochastic Processes]{Hidden Geometry of Bi-Directional Grid Constrained Stochastic Processes}

\author{Aldo Taranto, Shahjahan Khan} %<---Phase 1
%\author{Aldo Taranto, Ron Addie, Shahjahan Khan} %<---Phase 2

%\email{Aldo.Taranto@}
%\address{Mathematics \& Statistics, University of Southern Queensland, Toowoomba, QLD 4350, Australia}

%\vspace*{1cm}

%\curraddr{Credience Corporation, www.credience.com}
%\curremail{Aldo.Taranto@credience.com}

%\vspace*{1cm}

%\author{Ron Addie}
%\email{Ron.Addie@} %<----------------ADD BACK IN

%\author{Shahjahan Khan}
%\email{Aldo.Taranto@ \quad \quad @usq.edu.au}
\email{Aldo.Taranto@usq.edu.au, Shahjahan.Khan@usq.edu.au} %<----Phase 1
%\email{Aldo.Taranto@, Ron.Addie@, Shahjahan.Khan@, $\quad$ @usq.edu.au} %<----Phase 2

%\urladdrname{www.usq.edu.au}

%    Information for third author
%\author{Ravinesh C. Deo}
%\email{Ravinesh.Deo@usq.edu.au}
%\vspace*{0.5cm}

\thanks{The first author was supported by an Australian Government Research Training Program (RTP) Scholarship.}
\thanks{We would like to thank A/Prof. Ron Addie of University of Southern Queensland for advice on refining this paper.}

\address{\text{ }\\
School of Sciences\\
University of Southern Queensland\\
Toowoomba, QLD 4350, Australia}
%\thanks{Support information for the third author.}

%    General info
\date{\today}



\begin{abstract}
Bi-Directional Grid Constrained (BGC) stochastic processes (BGCSP) are constrained It\^{o} diffusions with the property that the further they drift away from the origin, the more resistance to movement in that direction they undergo.
We investigate the underlying characteristics of the BGC parameter $\Psi (x, t)$ by examining its geometric properties.
The most appropriate convex form for $\Psi$, i.e. the parabolic cylinder is identified after extensive simulation of various possible forms.
The formula for the resulting hidden reflective barrier(s) is determined by comparing it with the simpler Ornstein-Uhlenbeck process (OUP).
Applications of BGCSP arise when a series of semipermeable barriers are present, such as regulating interest rates and chemical reactions under concentration gradients, which gives rise to two hidden reflective barriers.
\end{abstract}

\maketitle %This makes the title page

%%%%%%%%%%%%%%%%%%%%%%%%%%%%%%%%%%%%%%%%%%%%%%%%%%%%%%%%%%%%%%%
%\newpage
%%%%%%%%%%%%%%%%%%\tableofcontents

%\vfill
%\vspace*{0.5cm}

%\bigskip
%%%%%%%%%%%%%%%%%%%%%%%%%%%%%%%%%%%%%%%%%%%%%%%%%%%%%%%%%%%%%%%
%\newpage
\begin{comment}
\section*{Notation}

\begin{table}[H]
   \begin{center}
\begin{tabular}{  l  l  }
\hline
 \textbf{Term} & \textbf{Description}   \\ \hline
BGC                & Bi-Directional Grid Constrained\\
BGCSP                & BGC Stochastic Process\\
$T$                 & Time\\
$t$                 & Index for time in $[0, T]$\\
$X_t$              & Stochastic process over time $t$\\
$W_t$              & Wiener process over time $t$ \\
$f(X_t, t)$, $\mu$             & Drift coefficient of $X$ over $t$\\
$g(X_t, t)$, $\sigma$         & Diffusion coefficient of $X$ over $t$\\
$\Psi(x, t)$         & BGC coefficient of $x$ over $t$\\
$\mathfrak{B}_L$  &  Hidden reflective lower barrier\\
$\mathfrak{B}_U$  &  Hidden reflective upper barrier\\
$| x |$            & Absolute value of $x$\\
%$X \perp Y$     & $X$ and $Y$ are perpendicular\\
%$\langle \bigcdot | \bigcdot \rangle$ & Hermitian inner product operator\\
%$\mathfrak{Re}(\bigcdot)$ & Real component of a complex number\\
%$\mathfrak{Im}(\bigcdot)$ & Imaginary component of a complex number\\
\hline
\end{tabular}
   \label{Tab:Notation}
   \end{center}
\end{table}
\end{comment}

\bigskip
%%%%%%%%%%%%%%%%%%%%%%%%%%%%%%%%%%%%%%%%%%%%%%%%%%%%%%%%%%%%%%%
\section{Introduction}

%\citep{TarantoKhan2020_1}
%\cite{TarantoKhan2020_1}

\noindent
In \citet{TarantoKhan2020_1}, the concept of Bi-Directional Grid Constrained (BGC) stochastic processes (BGCSP) was described as a general It\^{o} diffusion in which the further it drifts away from the origin, the more constrained the It\^{o} diffusion(s) becomes.
We note that for an arbitrary stochastic function $h$, that the following notations are equivalent,

\[
h(X(t)) = h(X, t) = h(X_t) = h(X_t, t)
\]

\bigskip \noindent
and that the last of these is adopted here.
We will also interchange $\Psi(X_t, t)$, $\Psi(x, t)$ and $\Psi(x)$ depending on the specific context.
The stochastic differential equation (SDE) of BGC stochastic processes was defined as follows.

\bigskip \noindent
\begin{definition}\textbf{(Definition I of BGC Stochastic Processes)}. For a complete filtered probability space $(\Omega, \mathcal{F}, \{ \mathcal{F} \}_{t \geq 0}, \mathbb{P})$ and a BGC function $\Psi (x) : \mathbb{R} \rightarrow \mathbb{R}$, $\forall x \in \mathbb{R}$, then the corresponding BGC It\^{o} diffusion is expressed as,

\bigskip \noindent
\begin{equation}\label{Eq:BGC1}
      dX_t  =  \Big( f(X_t, t)  \underbrace{-  \sgn[X_t, t] \Psi  (X_t, t) }_{\textbf{BGC}} \Big) \, dt + g(X_t, t) \, dW_t, 
\end{equation}

\bigskip \noindent
where $\sgn[x]$ is the sign function defined in the usual sense, $f(X_t, t)$ is the drift term, $\Psi (x, t)$ is the constraining term, $g(X_t, t)$ is the diffusion term and $f(X_t, t)$, $\Psi (x, t)$, $g(X_t, t)$ are convex functions. \hfill    $\blacksquare$
\end{definition}

\begin{figure}[H]
   \centering
\includegraphics[scale=0.26]{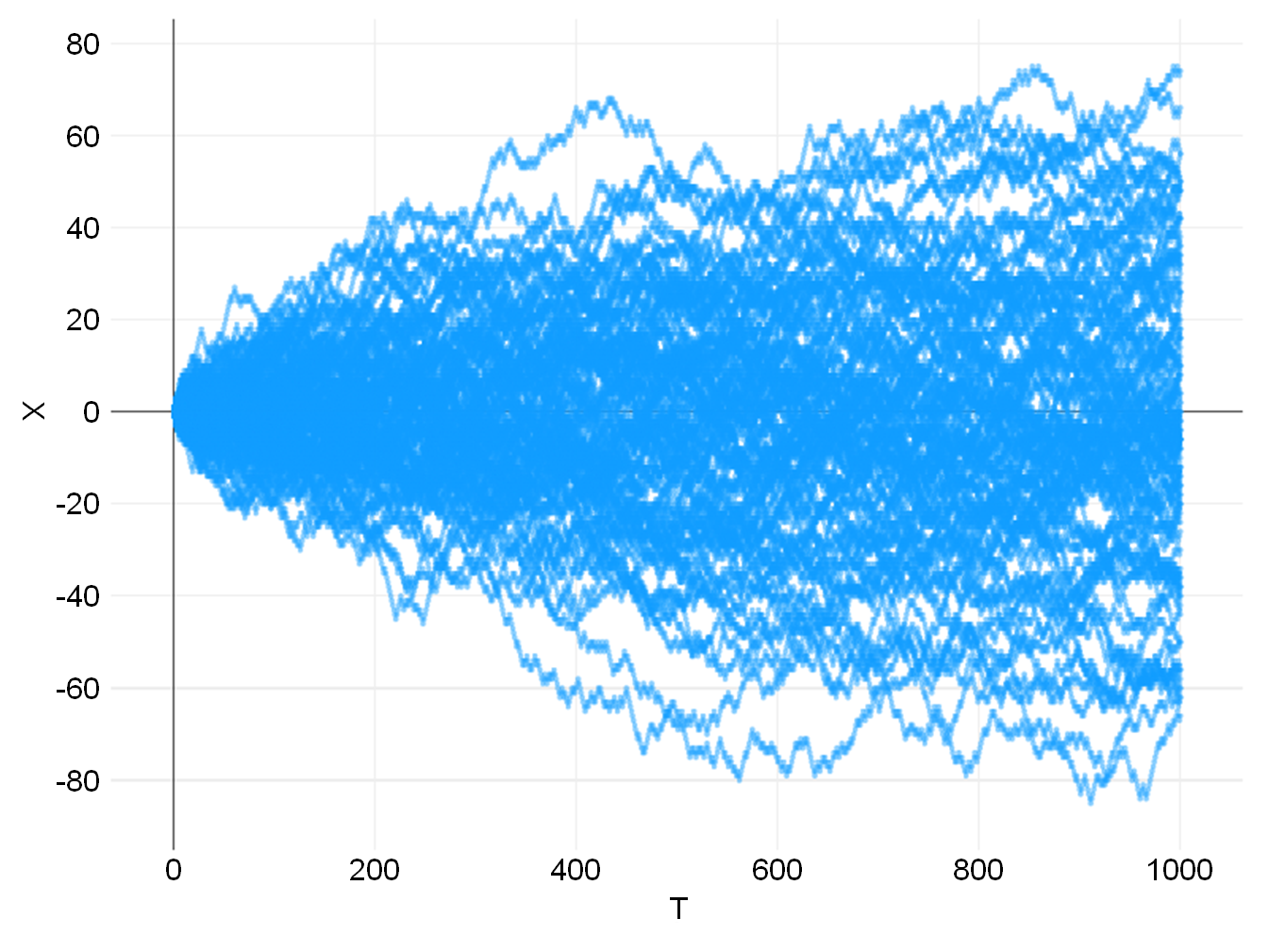}
\includegraphics[scale=0.26]{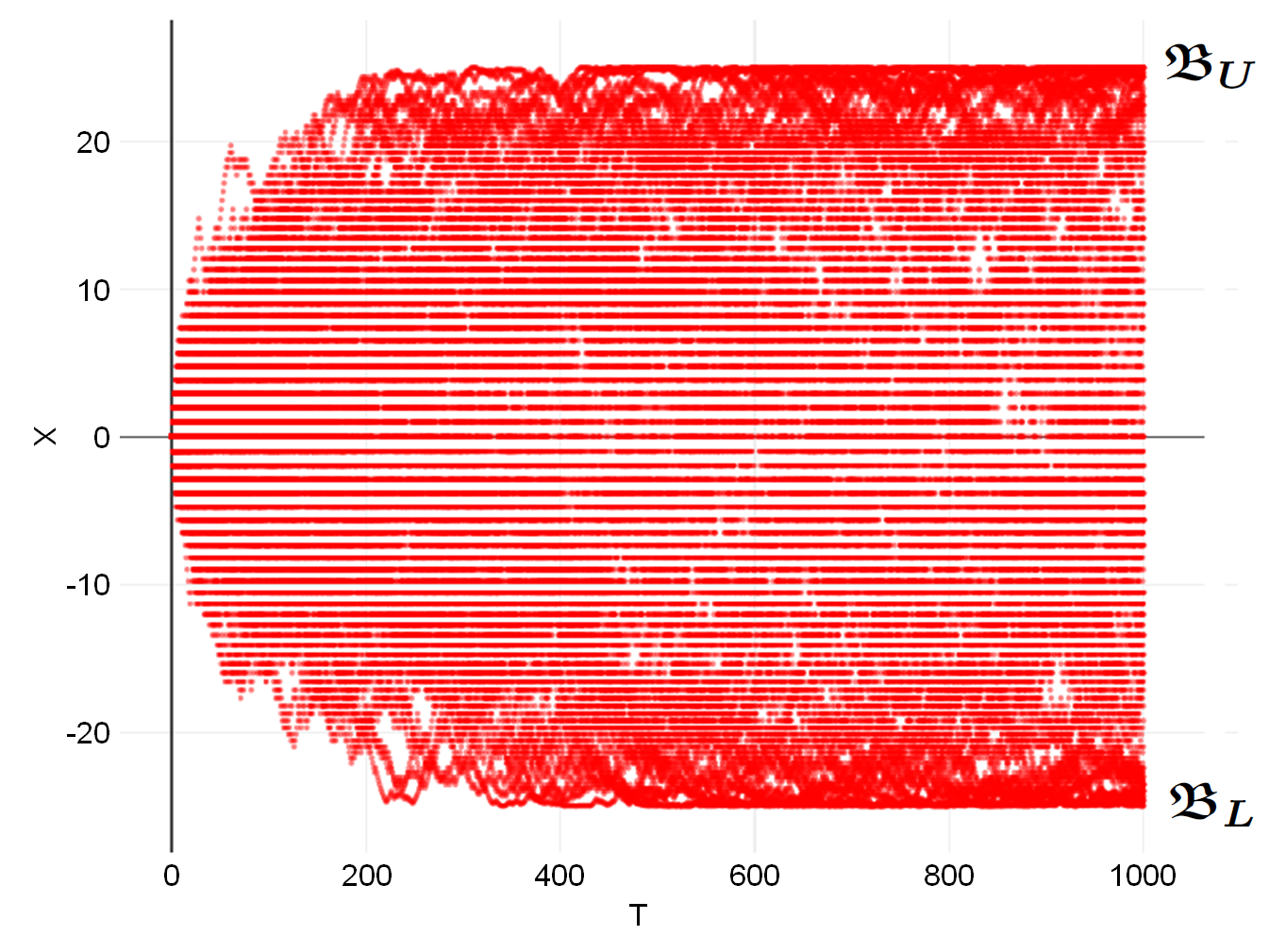}\\
   \textbf{\footnotesize \noindent{
(a). $\boldsymbol{\mu = 0}$, $\boldsymbol{\sigma = 1}$, without BGC  \quad \quad \quad \quad \quad (b). $\boldsymbol{\mu = 0}$, $\boldsymbol{\sigma = 1}$, with BGC \quad \quad \\}}
   \caption{It\^{o} Diffusions With \& Without BGC}
   \label{Fig:Grid_00}
\flushleft
   \textbf{\footnotesize The zero drift in (a) is constrained in (b) the more it deviates from the origin, causing the hidden reflective upper barrier $\boldsymbol{\mathfrak{B}_U}$ and hidden reflective lower barrier $\boldsymbol{\mathfrak{B}_L}$ to emerge, together with horizontal bands to form due to the discretization effect of BGC.
}
\end{figure}

\bigskip \noindent
%We notice that our original use of $\lambda (x)$ has been replaced in this paper with $\Psi (x)$, as $\lambda (x)$ is reserved for usage in the exponential distribution.
To visualize the impact of BGC stochastic processes, 1000 It\^{o} diffusions were simulated both with and without BGC, with unit diffusion coefficient $\sigma (X_t, t)$ for negative, zero and positive drift $\mu (X_t, t)$ coefficients.
Figure \ref{Fig:Grid_00} shows this when $\mu (X_t, t)=0$ and $\sigma (X_t, t)=1$, so that one can see the hidden upper barrier $\mathfrak{B}_U$ and hidden lower barrier $\mathfrak{B}_L$ emerge more clearly than is the case when using other coefficients.

\bigskip \noindent
\begin{remark}
The drift $f(X_t, t)$ and diffusion (or volatility) $g(X_t, t)$ terms reflect the instantaneous mean and standard deviation respectively.
It must also be noted from Figure \ref{Fig:Grid_00} that even when a generalized It\^{o} diffusion is reduced to a single Wiener process by setting $f(X_t, t) = 0$ and $g(X_t, t) = 1$, then BGC still impacts the the stochastic process.
\end{remark}

\bigskip \noindent
It is from these observations that an alternative definition to (\ref{Eq:BGC1}) can be stated as follows.

\bigskip \noindent
\begin{definition}\textbf{(Definition II of BGC Stochastic Processes)}. For a complete filtered probability space $(\Omega, \mathcal{F}, \{ \mathcal{F} \}_{t \geq 0}, \mathbb{P})$ and a BGC function $\Psi (x) : \mathbb{R} \rightarrow \mathbb{R}$, $\forall x \in \mathbb{R}$, then the corresponding BGC It\^{o} diffusion can be expressed as,

\bigskip \noindent
\begin{equation}\label{Eq:BGC2}
 dX_t  =   f(X_t,t)  \, dt  +  \Big( g(X_t, t) \underbrace{-  \sgn[X_t, t] \Psi  (X_t, t) }_{\textbf{BGC}} \Big) \, dW_t,
\end{equation}

\bigskip \noindent
where $\sgn[x]$ is defined in the usual sense, $f(X_t, t)$ is the drift term, $g(X_t, t)$ is the diffusion term, $\Psi (x, t)$ is the constraining term and $f(X_t, t)$, $g(X_t, t)$, $\Psi (x, t)$ are convex functions. \hfill    $\blacksquare$
\end{definition}

\bigskip \noindent
\begin{remark}
In \citet{TarantoKhan2020_1}, only 2 decimal places were used and some readers may argue that this is a low level of precision for simulation results to be robust.
To show that the discretization effect present in BGC stochastic processes is not due to such rounding errors, all our simulations were rerun to ten times more precision (i.e. to 20 decimal places) and it was found that the discretization or banding effect of BGC was still present, so it is a real phenomenon.
More shall be discussed about this in the Results and Discussion section.
\end{remark}

\bigskip \noindent
This paper will answer two main objectives;

\begin{enumerate}[(1).]
\item What are the key properties of $\Psi (X_t, t)$ that are relevant in BGC stochastic processes?
\item What is the formula for the hidden reflective lower barrier $\mathfrak{B}_L$ and the hidden reflective upper barrier $\mathfrak{B}_U$ in relation to $\Psi (X_t, t)$?
\end{enumerate}

\bigskip \noindent
Before these objectives are addressed in the Methodology section, the relevant research is examined in the Literature Review section.

%%%%%%%%%%%%%%%%%%%%%%%%%%%%%%%%%%%%%%%%%%%%%
%%%%%%%%%%%%%%%%%%
\bigskip
%\newpage
\section{Literature Review}

\noindent
\textbf{Constraining Discrete Random Walks}. Constrained stochastic processes have been applied to game theory (\citet{Feller1968}) and conditional Markov chains of this type have also been applied to biology, branching processes (\citet{FerrariMartinezPicco1992}), molecular physics (\citet{NovikovFieremansJensenHelpern2011}), medicine (\citet{Bell1976}) and queuing theory (\citet{BohmGopal1991}) to name a few.
\citet{Weesakul1961} discussed the classical problem of random walks restricted between a reflecting and an absorbing barrier.
\citet{Lehner1963} studied 1-Dimensional random walks with a partially reflecting barrier using combinatorial methods.
\citet{Gupta1966} introduced the concept of a multiple function barrier (MFB) where a state can either absorb, reflect, let through (transmit) or hold for a moment along with their corresponding probabilities.
\citet{DuaKhadilkarSen1976} found the bivariate generating functions of the probabilities of a random variable reaching a certain state under different conditions.
\citet{Percus1985} considered asymmetric random walks, with one or two boundaries, on a 1-Dimensional lattice.
\citet{ElShehawey2000} obtained absorption probabilities at the boundaries for  random walks between one or two partially absorbing boundaries, using conditional probabilities.

\bigskip \noindent
\textbf{Constraining Continuous Wiener Processes}. Dirichlet studied the first boundary value problem, for the Laplace equation, proving the uniqueness of the solution and this type of problem in the theory of partial differential equations (PDEs).
This was later named the Dirichlet problem after him (\citet{GowersBarrowGreenLeader2008}).
Problems expressed within this framework were studied as early as 1840 by C.F. Gauss, and then by \citet{Dirichlet1850}.
\citet{Kurtz1991} formulated a means for constraining Markov processes.
\citet{Lepingle2009} expanded upon previous research on barriers, which included boundary behavior of constrained Wiener processes between reflecting and repellent barriers.
\citet{MajumdarRandonFurlingKearneyYor2008} derived the time taken to reach the maximum for a variety of constrained Wiener processes.
\citet{OrmeciDaiVate2008} examined the constraining of Wiener processes via impulse control.
\citet{BudhirajaDupuis1999} added necessary and sufficient conditions for the stability of such constrained processes.
The same authors studied large deviations for various metrics of reflecting Wiener processes under constraining (\citet{BudhirajaDupuis2003}).
\citet{KharroubiMaPhamZhang2010} constrained the jumps of Backward SDEs.

\bigskip \noindent
Whilst BGC stochastic processes are relatively new, they do have applications in many areas, most prominent being in mathematical finance, investment algorithms and quantitative trading (Taranto and Khan, \cite{TarantoKhan2020_4}, \cite{TarantoKhan2020_5}, \cite{TarantoKhan2020_6}, \cite{TarantoKhan2020_7}).

\bigskip \noindent
We are now in a position to examine the geometry of the random variable $X$ as constrained by $\Psi(X_t, t)$ and its ramifications for BGCSP.

%%%%%%%%%%%%%%%%%%%%%%%%%%%%%%%%%%%%%%%%%%%%%
%%%%%%%%%%%%%%%%%%
\bigskip
%\newpage
\section{Methodology}

\subsection{Convexity of BGC}

\noindent
From (\ref{Eq:BGC2}), we know that $\Psi (X_t, t)$ needs to be a convex function and specifically, {\it centered about the origin}.
This is to ensure that the constraining applies increasing monotonic resistance to the It\^{o} diffusion in both directions (i.e. bi-directionally).
For example, $\Psi (x, t) = e^{t}$ would not be sufficient because whilst $e^{t}$ is convex (as shown by having a line bisect any two points on its curve) $e^{t}$ does not increase as $t \downarrow -\infty$ as it does when $t \uparrow +\infty$. 
This gives rise to the need for the following classification of convexity.

\bigskip \noindent
\begin{definition}
\textbf{(Types of Convexity)}.
If $\Psi (x): \mathbb{R} \rightarrow \mathbb{R}$ and $\Psi (x) \in C^2$, then we can characterize (\citet{Zalinescu2002}, \citet{BauschkeCombettes2011}) its convexity as follows,

\bigskip \noindent
\begin{enumerate}
\item $\Psi (x)$ is {\it \textbf{convex}} if and only if $\Psi'' (x) \geq 0$, $\forall x \in \mathbb{R}$.
\item $\Psi (x)$ is {\it \textbf{strictly convex}} if and only if $\Psi'' (x) > 0$, $\forall x \in \mathbb{R}$.
\item $\Psi (x)$ is {\it \textbf{strongly convex}} if and only if $\Psi'' (x) \geq m > 0$, $\forall x \in \mathbb{R}$.  \hfill    $\blacksquare$
\end{enumerate}
\end{definition}

\bigskip \noindent
We will require a new type of (subset) convexity for BGCSP.

\bigskip \noindent
\begin{definition}
\textbf{(Bi-Directional Convexity)}.
If $\Psi (x): \mathbb{R} \rightarrow \mathbb{R}$ and $\Psi (x) \in C^2$, then we can characterize its convexity as follows,

\smallskip
$\Psi (x)$ is {\it \textbf{bi-directionally convex}} if and only if $\Psi'' (x) > 0$, $\Psi(x)= \Psi(-x)$, $\forall x \in \mathbb{R}$.
 \hfill    $\blacksquare$
\end{definition}

\begin{figure}[H]
   \centering
\includegraphics[width=\linewidth]{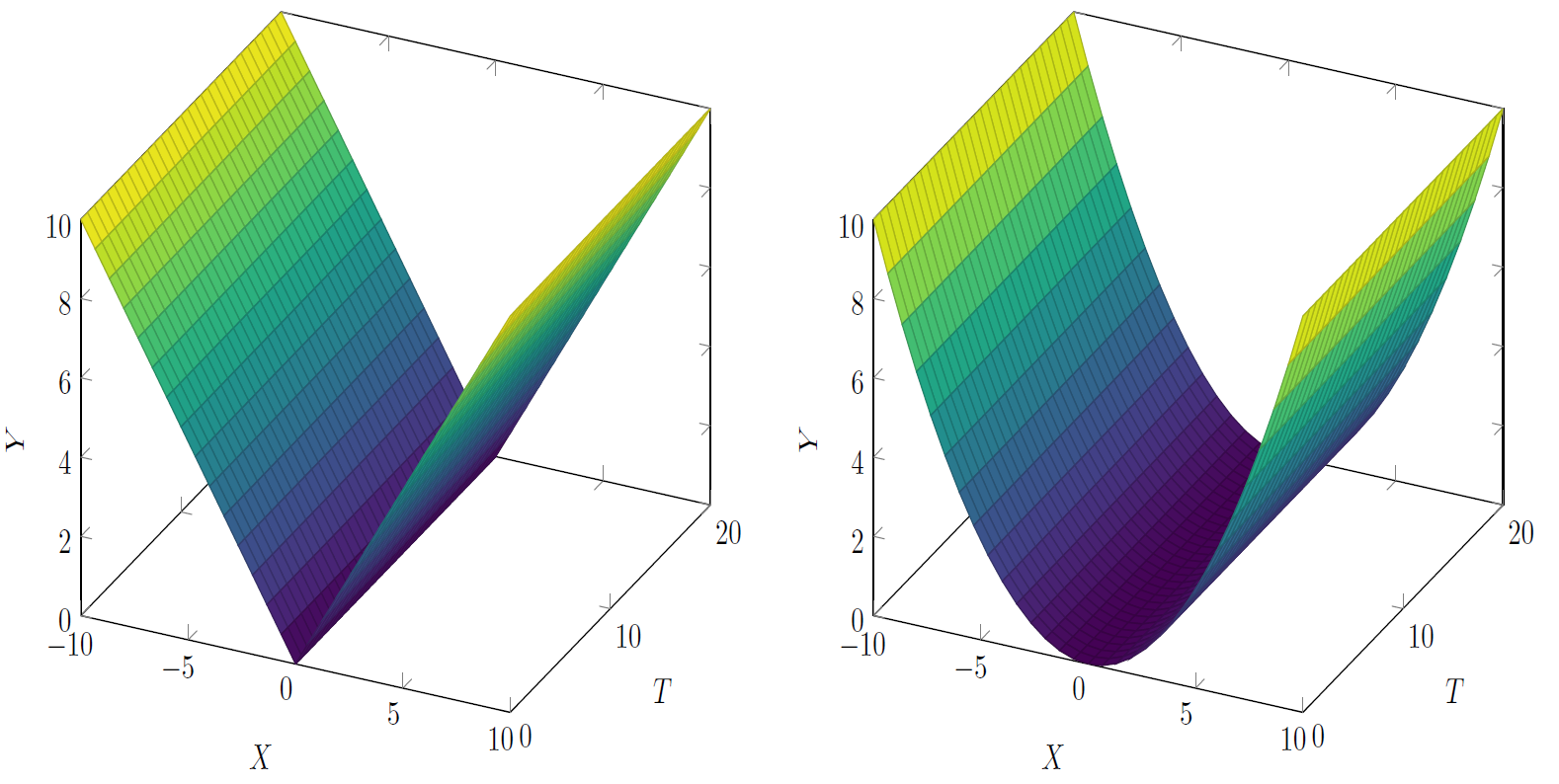}
   \textbf{\footnotesize \noindent{
(a). $\boldsymbol{y = |x|}$ \quad \quad \quad \quad \quad \quad \quad \quad \quad \quad \quad \quad \quad (b). $\boldsymbol{y = x^2 / \omega}$, $\boldsymbol{\omega =10}$\\}}
\includegraphics[width=\linewidth]{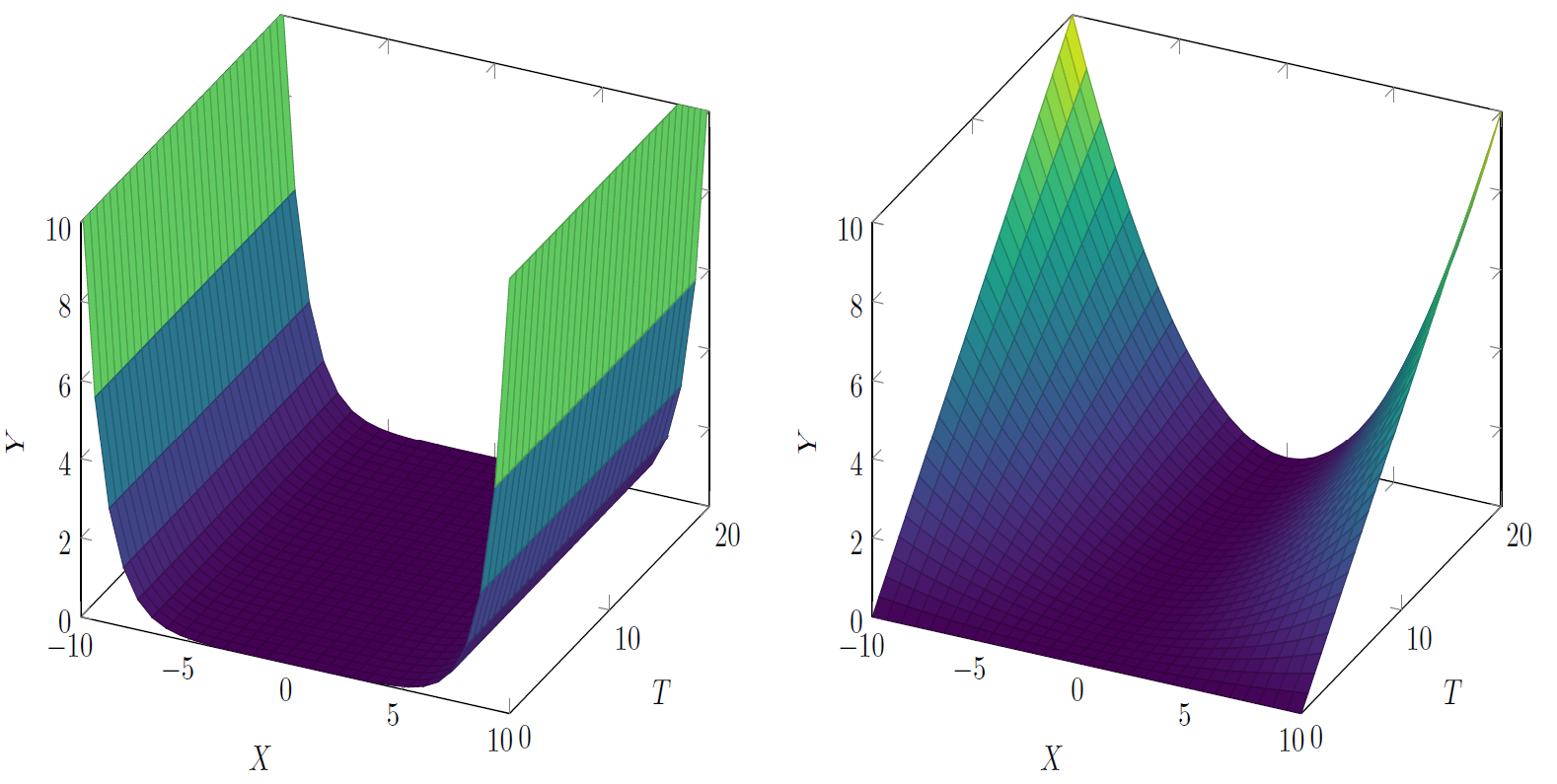}
   \textbf{\footnotesize \noindent{
(c). $\boldsymbol{y = (e^x + e^{-x})/ \omega}$, $\boldsymbol{\omega=2000}$ \quad \quad \quad \quad \quad \quad \quad    (d). $\boldsymbol{y = x^2 t/\omega}$, $\boldsymbol{\omega=200}$\\}}
%=============================================
   \caption{Use Cases I - IV: Surface Plots of Potential Convex BGC Function $\Psi (x, t)$}
   \label{Fig:Grid_02}
\flushleft
   \textbf{\footnotesize
(a). \textbf{USE CASE I}: Wedge.\\
(b). \textbf{USE CASE II}: Parabolic cylinder.\\
(c). \textbf{USE CASE III}: Double exponential cylinder.\\
(d). \textbf{USE CASE IV}: Hybrid of flat plane and parabolic cylinder.\\
}
\end{figure}

\begin{figure}[H]
   \centering
\includegraphics[width=\linewidth]{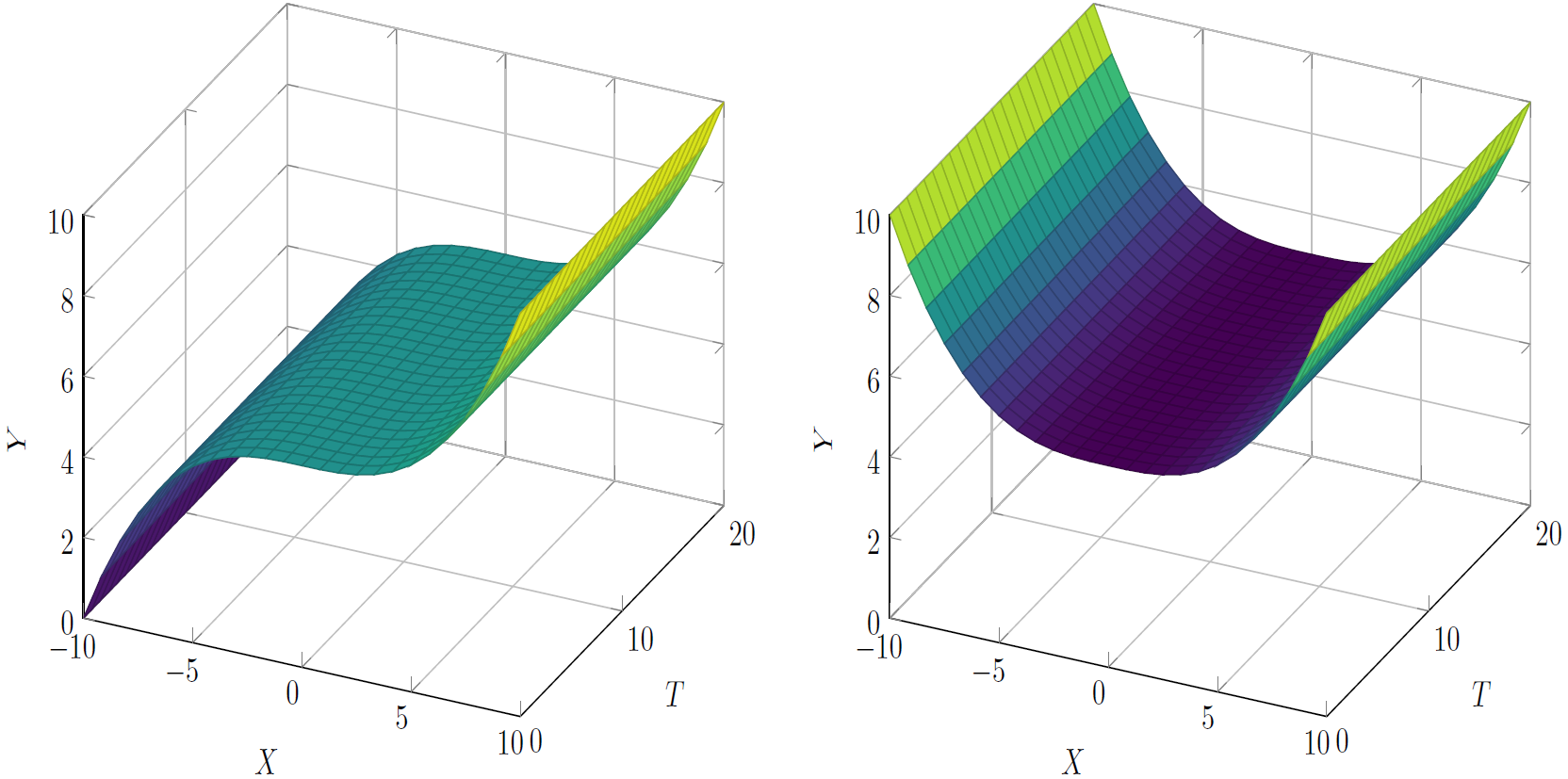}
   \textbf{\footnotesize \noindent{
(a). $\boldsymbol{y = x^3 / \omega_1 + \omega_2}$, $\boldsymbol{\omega_1 = 200}$, $\boldsymbol{\omega_2=5}$         \quad \quad  (b). $\boldsymbol{y = |x|^3 / \omega_1 + \omega_2}$, $\boldsymbol{\omega_1 = 200}$, $\boldsymbol{\omega_2=5}$}}
%=============================================
   \caption{Technique to Generate Convex BGC Polynomial Cylinders}
   \label{Fig:Grid_03}
   \textbf{\footnotesize
(a). Cubic cylinder, $\quad$ (b). \textbf{USE CASE V} - Spliced polynomial cylinder(s).\\
\flushleft
(a). The cubic is concave for negative $\boldsymbol{X}$ values and convex for positive $\boldsymbol{X}$ values.\\
\flushleft
(b). By forcing the negative values to be positive by the absolute value function, then the entire function becomes convex and is Bi-Directional of the polynomial cylinder variety (in this case a cubic cylinder).
It is as if the positive part of the cubic cylinder was spliced in to replace the negative part.\\
}
\end{figure}

\bigskip \noindent
To establish some use cases to explore the convex geometry of potential BGC functions for $\Psi (x, t)$, we plot their surfaces in Figures \ref{Fig:Grid_02} and \ref{Fig:Grid_03}, and deduce which specific type of convexity definition is required for BGCSP.

\bigskip \noindent
We examine five main tpes of bi-directionally convex functions in $\mathbb{R}^3$, as shown in Figures \ref{Fig:Grid_02} and \ref{Fig:Grid_03}.

\begin{enumerate}
\item Use Case I was not adopted because it constrains the It\^{o} diffusion evenly and uniformly but not rapidly enough with the unconstrained It\^{o} diffusions' iterated logarithm bounds growth rate.
This is to the point that the before and after BGC plots (see Figure \ref{Fig:PotentialConvexCandidatesforPsi}(a)) look very similar and do not constitute a practical and worthwhile BGC process.

\item Use Case II is definitely the ideal function for BGC and so we will dedicate much of the Methodology, Results and Discussion sections to the parabolic cylinder.

\item Use Case III examines how the standard convex $e^{x}$ will not suffice because $\Psi (x, t)$ needs to be `Bi-Directional', that is, it must be a mirror reflection about the origin $X = 0$ over all time.
Hence, $e^x + e^{-x}$ was used, yet due to the fast growing nature of the exponential function, it constrains the It\^{o} process too much for it to be a useful function for BGC (even when it is scaled down by a constant $\omega$ or many other possible variations of $e^{\pm x}$), as will be elaborated further in the Results and Discussion section Figure \ref{Fig:PotentialConvexCandidatesforPsi}(b).

\item Case IV was also presented here because it is a transition from no constraining to a gradual parabolic cylinder.
Such a surface was proposed for applications in which It\^{o} diffusions are not constrained so much initially and which become increasingly more constrained over time are required.
However, as will be detailed further in the Results and Discussion section Figure \ref{Fig:PotentialConvexCandidatesforPsi}(c), this did not produce the hidden barriers that bound an It\^{o} diffusion from both above and below.

\item Use Case V is the final example that is worthwhile discussing, as shown in Figure \ref{Fig:Grid_03}.
\end{enumerate}

\bigskip \noindent
From Figure \ref{Fig:Grid_03}, we note that $y = x^{2n}$, $\forall n \in \mathbb{N}$, will always result in a polynomial cylinder as $2n$ will always be an even exponent.
For odd exponents $y = x^{2n+1}$, $\forall n \in \mathbb{N}$, one can simply replace this with $y = |x|^{2n+1}$.
In general, $y=|x|^n$ will be a polynomial cylinder that will always be convex and `Bi-Directional', $\forall n \in \mathbb{N}$.
As will be elaborated in the Results and Discussion section in Figure \ref{Fig:PotentialConvexCandidatesforPsi}(d), the polynomial cylinder was not suitable as a BGC function for general unconstrained It\^{o} diffusions, but can be scaled to suit one's specific unconstrained It\^{o} diffusion.

%\FloatBarrier

%\FloatBarrier

\bigskip
\begin{remark}
It is clear by now that not any convex function can be appropriate for BGC.
An example of this would be $y = x^2 + t^2$, where it is clearly and bi-directionally convex but not constant or `cylindrical' over time and does not resemble any natural regime to constrain the stochastic processes uniformly over time.
It is thus clear now that BGC requires the bi-directionally convex definition and in particular, bi-directionally convex cylinders.
\end{remark}

\begin{figure}[H]
   \centering
\includegraphics[width=\linewidth]{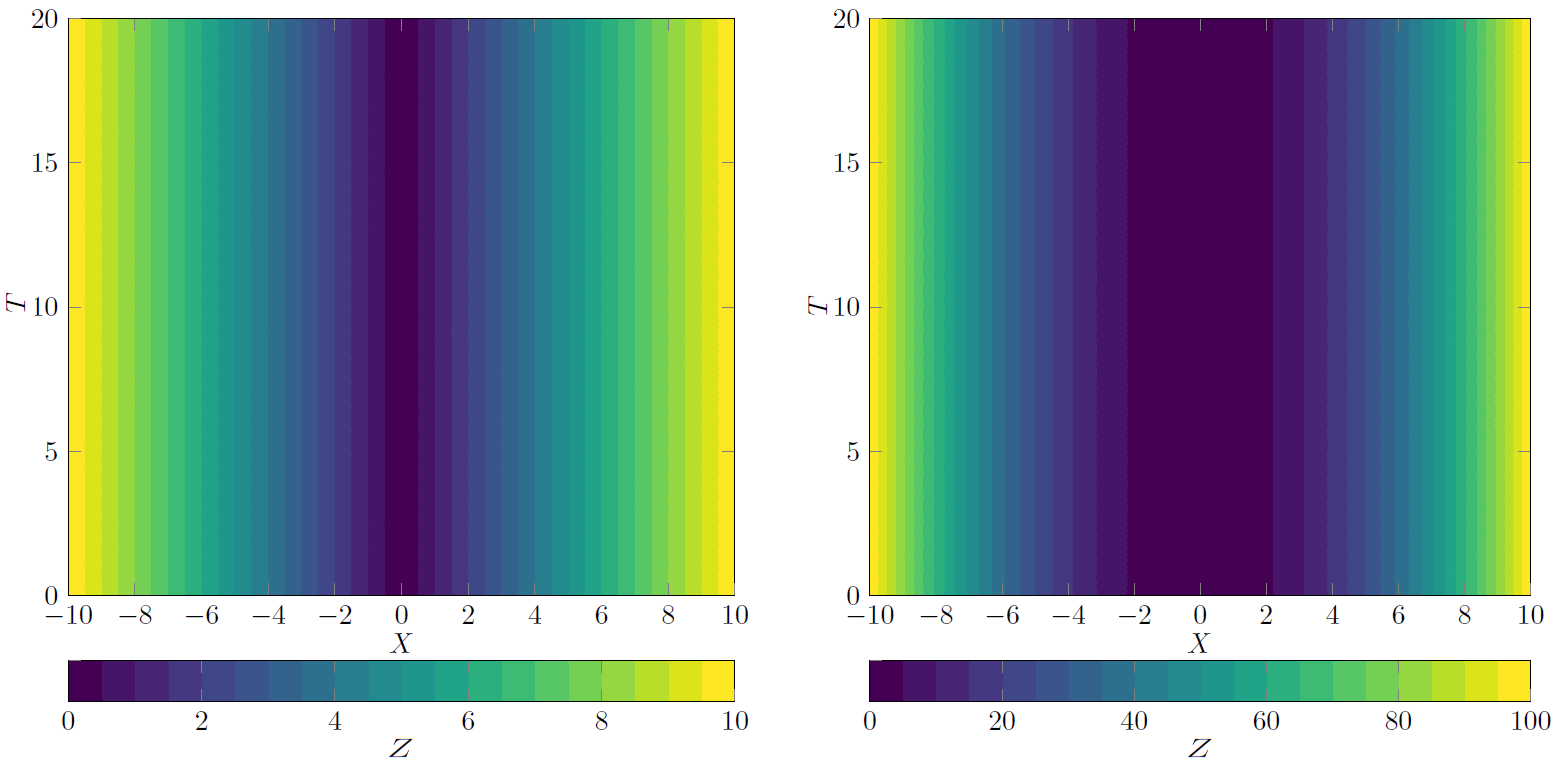}
   \textbf{\footnotesize \noindent{
(a). \quad \quad $\boldsymbol{y = |x|}$   \quad \quad \quad \quad \quad \quad \quad \quad \quad \quad (b). $\boldsymbol{y = x^2 /\omega, \text{ } \omega=10}$}}\\
\includegraphics[width=\linewidth]{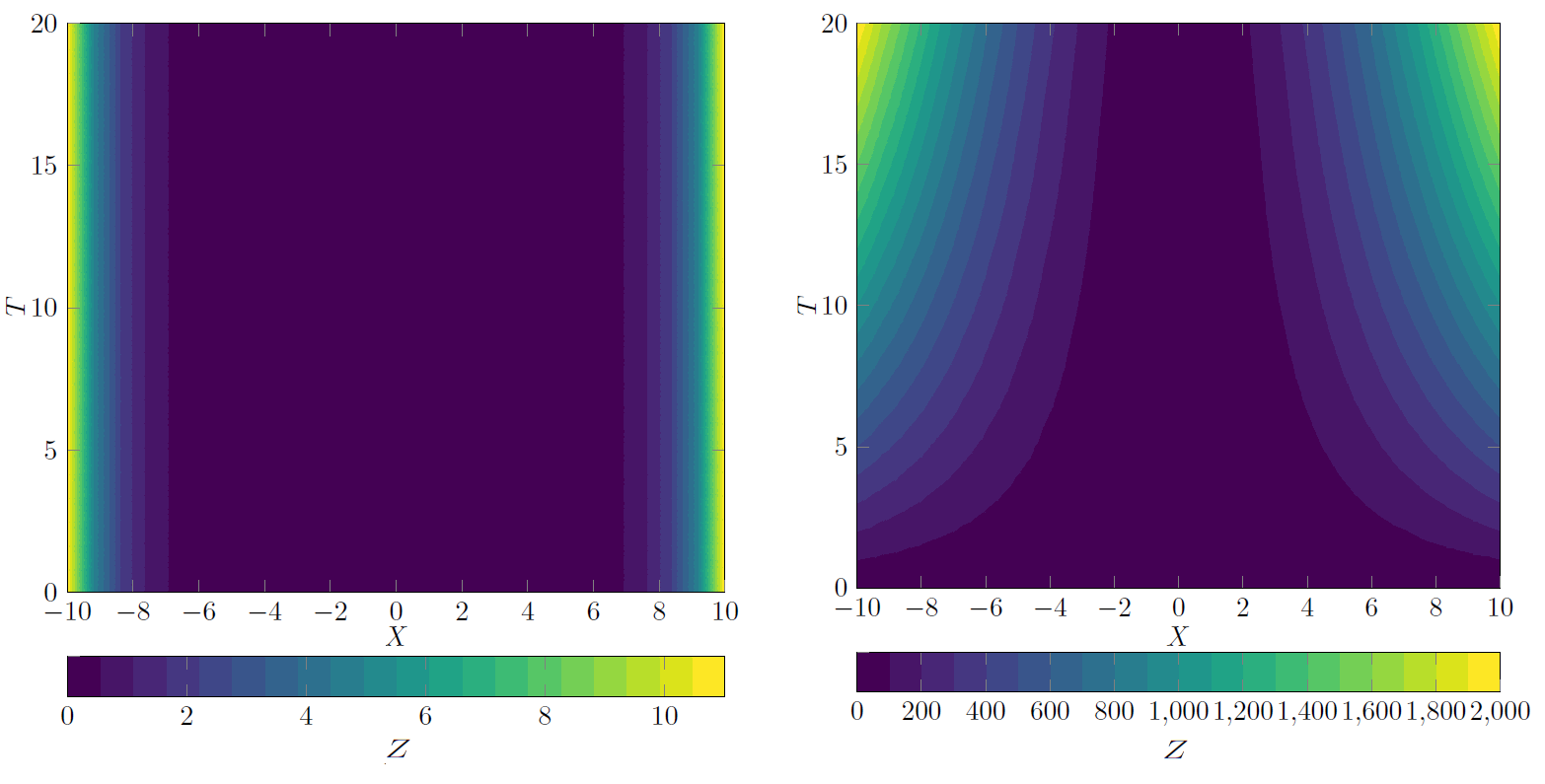}
   \textbf{\footnotesize \noindent{
(c). $\boldsymbol{y = (e^x + e^{-x})/\omega, \text{ } \omega=2000}$   \quad \quad \quad \quad \quad \quad \quad (d). $\boldsymbol{y = x^2 t/\omega, \text{ } \omega=200}$}}
%=============================================
   \caption{Contour Plots of Main Use Case Candidates for BGC}
   \label{Fig:Grid_04}
\flushleft
   \textbf{\footnotesize
The $\boldsymbol{Z}$ scale shows the 3D surface's height as it is mapped onto $\boldsymbol{\mathbb{R}^2}$ and $\boldsymbol{\omega}$ is a scaling function.\\
(a). \textbf{USE CASE I}: Contour of a Wedge.\\
(b). \textbf{USE CASE II}: Contour of a Parabolic Cylinder.\\
(c). \textbf{USE CASE III}: Contour of a Double Exponential Cylinder.\\
(d). \textbf{USE CASE IV}: Contour of a Hybrid of Flat Plane and Parabolic Cylinder.\\
}
\end{figure}

\begin{figure}[H]
   \centering
\includegraphics[width=\linewidth]{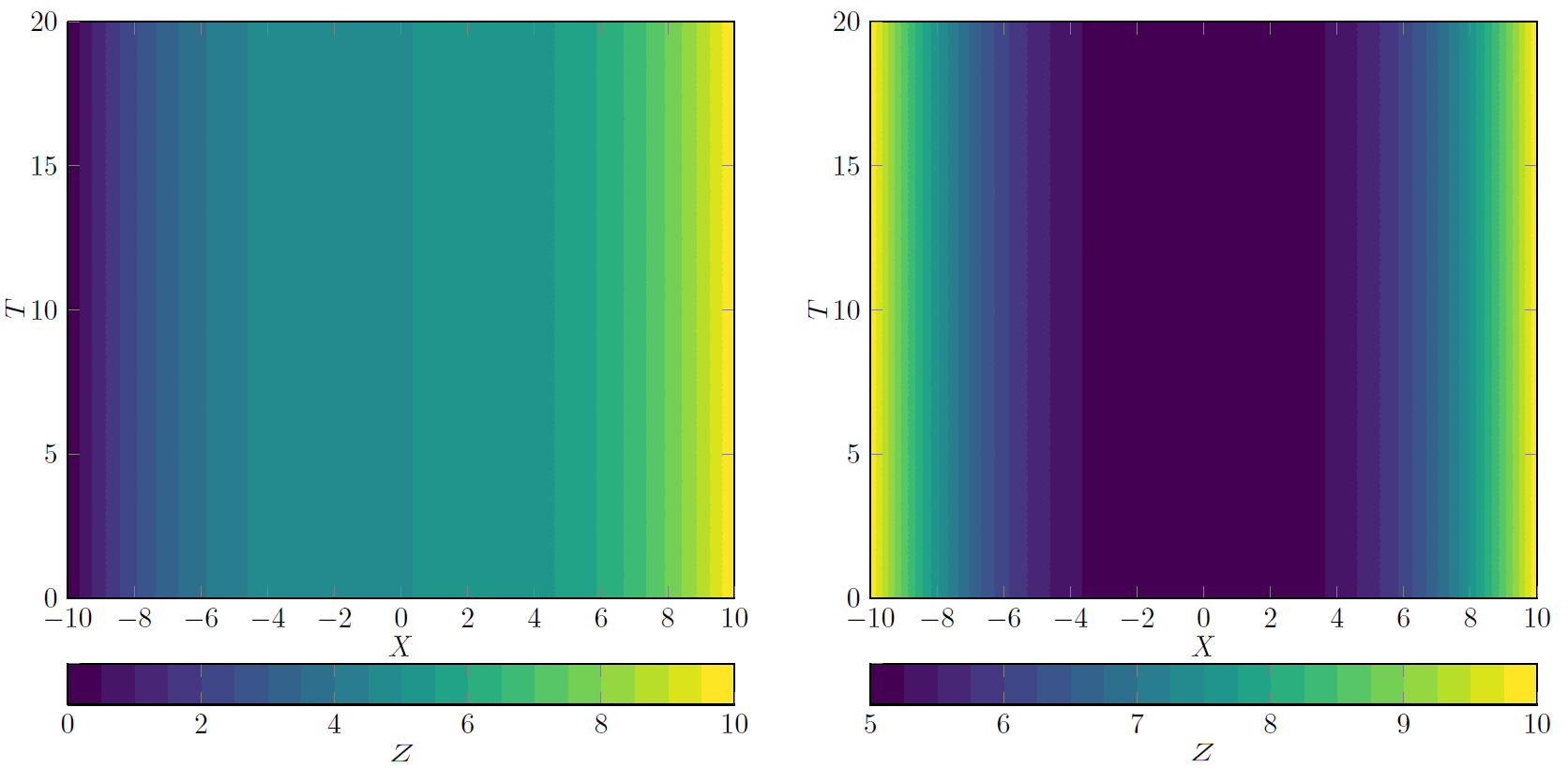}
   \textbf{\footnotesize \noindent{
(a). $\boldsymbol{y = x^3 / \omega_1 + \omega_2, \text{ } \omega_1 = 200, \text{ } \omega_2=5}$   \quad \quad  (b). $\boldsymbol{y = |x|^3 / \omega_1 + \omega_2, \text{ } \omega_1 = 200, \text{ } \omega_2=5}$}}
%=============================================
   \caption{Contour Plots as $y=x^3$ is made Bi-Directional}
   \label{Fig:Grid_05}
\flushleft
   \textbf{\footnotesize
The $\boldsymbol{Z}$ scale shows the 3D surface's height as it is mapped onto $\boldsymbol{\mathbb{R}^2}$.\\
(a). \boldsymbol{$y=x^3$}: We can see how negative values of $\boldsymbol{X}$ lead to a concave contour, whereas positive values of $\boldsymbol{X}$ lead to a convex contour, due to the odd power of $\boldsymbol{x^3}$. \\
(b). \textbf{USE CASE V}: By taking the absolute value of $\boldsymbol{x}$ before passing it through $\boldsymbol{x^3}$, we see that the contour plot is convex in both directions of $\boldsymbol{X}$, i.e. is Bi-Directional.\\
}
\end{figure}

\bigskip \noindent
Having explored the nature of $X$ as determined by $\Psi(x, t)$ which lies in $\mathbb{R}^3$, we notice that our It\^{o} process in $X$ is a 1-Dimensional stochastic process in $\mathbb{R}$, which when it propagates over time, it does so in $\mathbb{R}^2$.
The way we can see how the 3-D $\Psi(x, t)$ constrains the 1-D It\^{o} process in 2-D is via the projection of $\Psi(x, t)$ onto the $\mathbb{R}^2$ plane is via contour plots, as shown in Figures \ref{Fig:Grid_04} and \ref{Fig:Grid_05}.

\bigskip \noindent
From Figure \ref{Fig:Grid_04}, it is much clearer to see how $\Psi(x, t)$ constrains the It\^{o} process as $X_t$ propagates over time, where the lighter the colour, then the greater the resistance and hence the greater the constraining impact due to BGC.
This is also shown in the contour plots of Figure \ref{Fig:Grid_05}.

\bigskip \noindent
From these contour plots, we see how the convexity forms a series of decreasing semipermeable barriers (i.e. increasing reflection) on the It\^{o} process.
We now examine the effect that this has on the actual hidden reflective barriers, $\mathfrak{B}_L$ and $\mathfrak{B}_U$.

\bigskip
\begin{remark}
Note that an alternative to the 3-D surface inducing the 2-D contours, is the 3-D surface inducing the 3-D and 2-D vector fields, as shown in Figure \ref{Fig:VectorFieldInducedbyBGCConstraints}.

\begin{figure}[H]
   \centering
 \includegraphics[width=\linewidth]{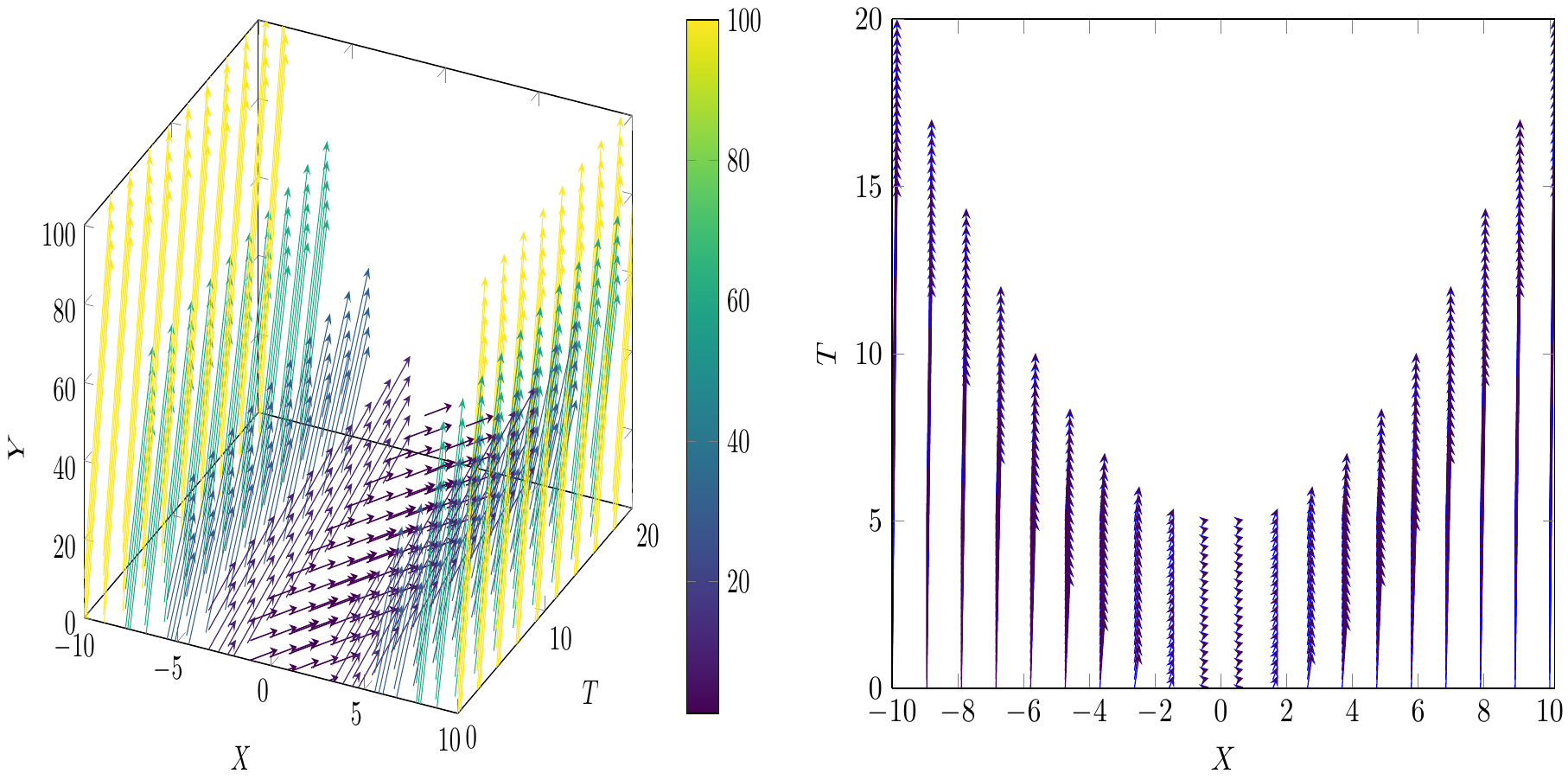}\\
   \textbf{\footnotesize \noindent{
(a). BGC Vector Field in $\boldsymbol{\mathbb{R}^3}$, $\quad \quad \quad$  (b). BGC Vector Field in $\boldsymbol{\mathbb{R}^2}$.
}}
   \caption{Vector Fields Induced by BGC Constraints}
   \label{Fig:VectorFieldInducedbyBGCConstraints}
\flushleft
   \textbf{\footnotesize
Just as in Figure \ref{Fig:Grid_04}(b) for the 2-D contour map induced by the 3-D parabolic cylinder surface, here the 3-D vector field in (a) also induces a constraining force on the 2-D vector plot in (b).
}
\end{figure}
%\FloatBarrier

\bigskip \noindent
From Figure \ref{Fig:VectorFieldInducedbyBGCConstraints}, we see that the constraining on $\mathbb{R}^2$ (specifically $\mathbb{R} \times \mathbb{R}_{+}$) can also be induced by the parabolic cylinder of $\Psi (X_t, t)$ and its associated vector field.
As the It\^{o} process propagates through the vector field, the greater the vector magnitudes, then the greater the resistance force of reflection back to the origin.

\bigskip \noindent
This novel concept has been researched recently, but in reverse by \citet{SimpsonKuske2018}, by modelling a constant variable into a random vector field to induce a stochastic process.
Specifically, they show how a Flippov system near a switching manifold (due to the meeting of vector fields) attracts orbits or constant variables in the abscence of randomness to create stochastic flow within the field.
  \hfill    $\blacksquare$
\end{remark}

%--------------------------------------------------------------------------------------------------------
\bigskip \noindent
\subsection{Hidden Barriers of BGC Stochastic Processes}

\begin{figure}[H]
   \centering
\includegraphics[width=\linewidth]{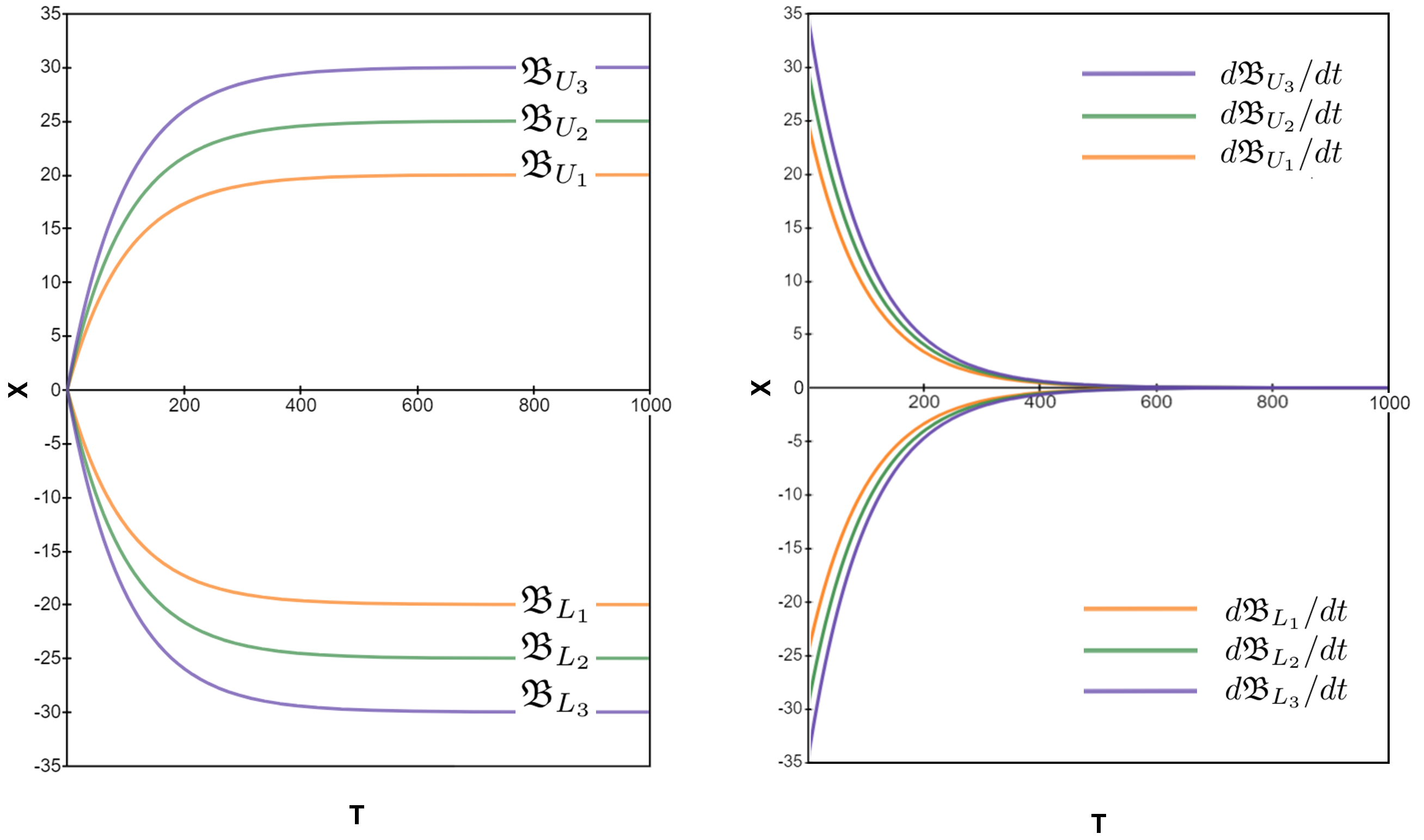}
   \textbf{\footnotesize \flushleft
$\quad \quad \quad \quad$ (a). $\quad \quad \quad \quad \quad \quad \quad \quad \quad \quad \quad \quad \quad \quad \quad \quad \quad \quad \quad$   (b).
}
%=============================================
   \caption{BGC Hidden Reflective Barriers Identified}
   \label{Fig:BGCHiddenReflectiveBarriersIdentified}
   \textbf{\footnotesize
Orange: $\boldsymbol{A= \pm 25}$, $\quad$ Green: $\boldsymbol{A= \pm 30}$, $\quad$ Purple: $\boldsymbol{A= \pm 35}$.\\}
\flushleft
   \textbf{\footnotesize
(a). $\boldsymbol{\mathfrak{B}_L = -A(1-e^{-\theta T})}$, $\boldsymbol{\mathfrak{B}_U = A(1-e^{-\theta T})}$, $\boldsymbol{\theta=0.01}$.\\
(b). The rate of change of the barriers, i.e. $\boldsymbol{d \mathfrak{B}_L / dt}$ and $\boldsymbol{d \mathfrak{B}_U / dt}$.
}
\end{figure}

\noindent
Whilst Figure \ref{Fig:Grid_00}(b) shows the detailed nature of the hidden reflective BGC varies, we only plot $\mathfrak{B}_L$ and $\mathfrak{B}_U$ to better help derive the formulation of the barriers, as shown in Figure \ref{Fig:BGCHiddenReflectiveBarriersIdentified}.

\bigskip \noindent
From Figure \ref{Fig:BGCHiddenReflectiveBarriersIdentified}, we see that the BGC hidden reflective barriers are regulated to one's desired distance from the origin by altering the $A$ parameter and are regulated in their climb rate from the origin to the barrier by altering the $\theta$ parameter.

%\FloatBarrier

\bigskip 
\begin{remark}
Note that these barriers are not the traditional constant reflective barriers such as $X= a$ or $X=b$ because the BGC It\^{o} diffusions are bounded within these barriers even as they depart from the origin (i.e. the unconstrained It\^{o} diffusions exceed these barriers, even near the origin), hence the initial curvature in the barriers.
This is generalized in the following Thorem.
\end{remark}

%===================================
\bigskip \noindent
\begin{theorem}
\textbf{(Hidden Barriers of BGC Stochastic Processes)}. For a complete filtered probability space $(\Omega, \mathcal{F}, \{ \mathcal{F} \}_{t \geq 0}, \mathbb{P})$ and a BGC function $\Psi (x) : \mathbb{R} \rightarrow \mathbb{R}$, $\forall x \in \mathbb{R}$ and the corresponding BGC It\^{o} diffusion expressed as,

\bigskip \noindent
%\[
\begin{equation}\label{Eq:BGCSP1}
      dX_t  =  \Big( \overbrace{ f(X_t, t)  \underbrace{-  \sgn[X_t, t] \Psi  (X_t, t) }_{\textbf{BGC}} }^{\boldsymbol{\mu (X_t, t)}} \Big) \, dt + \overbrace{g(X_t, t)}^{\boldsymbol{\sigma (X_t, t)}}  \, dW_t, 
\end{equation}
%\]

\bigskip \noindent
for $t \in [0,T]$, where $\sgn[x]$ is the sign function defined in the usual sense, $f(X_t, t)$ is a drift term, $\Psi (x, t)$ is the BGC term, $g(X_t, t)$ is the diffusion term and $f(X_t, t)$, $\Psi (x, t)$, $g(X_t, t)$ are bi-directionally convex functions.
Then the hidden lower barrier $\mathfrak{B}_L$ and hidden upper barrier $\mathfrak{B}_U$ are given by,

\begin{equation}\label{barriers}
\mathfrak{B}_L \geq -A(1-e^{-\theta T}), \quad \mathfrak{B}_U \leq A(1-e^{-\theta T}),
\end{equation}

\bigskip \noindent
where $A$,$\theta \in \mathbb{R}$ are constants, $A$ is the distance from the origin to the barrier(s) and $\theta$ is the rate of growth towards the barrier(s).
\end{theorem}

\bigskip
\begin{proof}
To some readers, (\ref{barriers}) is obvious just by looking at Figure \ref{Fig:BGCHiddenReflectiveBarriersIdentified}(a) simply because the function must asymptote horizontally (exponentially) towards the barrier(s).
However, this does not constitute a proof because we have a stochastic (and not a deterministic) process.
Assume for a moment that the above BGC SDE is a simpler object in which $\mu (X_t, t)$ and $\sigma (X_t, t)$ are constant, where the drift function $\mu(x): \mathbb{R} \rightarrow \mathbb{R}$ and the diffusion function $\sigma(x): \mathbb{R} \rightarrow \mathbb{R}$, $\forall x \in \mathbb{R}$, in the limit approach the typical constant expressions for the drift and diffusion coefficients, $\lim_{x \rightarrow \infty} \mu(x) \rightarrow \mu$, $\lim_{x \rightarrow \infty} \sigma(x) \rightarrow \sigma$.
We now start to see the resembelence of our BGC It\^{o} process with the simpler Ornstein–Uhlenbeck process (OUP),

\begin{equation}\label{Ornstein_Uhlenbeck_Process}
 dX_t = \overbrace{\kappa (\alpha - X_{t})}^{\boldsymbol{\mu (X_t, t)}} \, dt+\sigma \, dW_{t},
\end{equation}

\bigskip \noindent
where $\alpha$ is the long-term mean, $\kappa$ is the `attraction rate' or speed of mean-reversion, both of which are constants, as depicted in Figure \ref{Fig:BGCHiddenReflectiveBarriersRelationtoOrnsteinUhlenbeckProcess}.

\begin{figure}[H]
   \centering
\includegraphics[width=\linewidth]{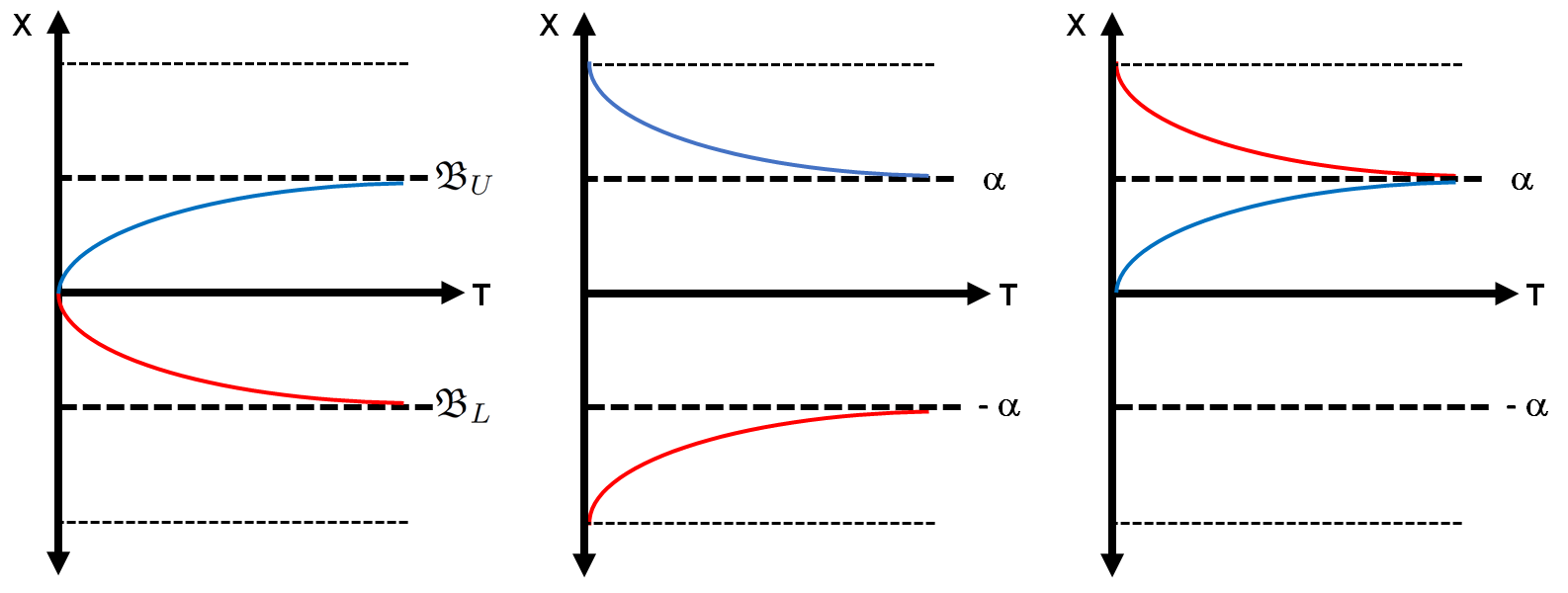}
   \textbf{\footnotesize \noindent{ \flushleft
(a). $\quad \quad \quad \quad \quad \quad \quad \quad \quad \quad$ (b). $\quad \quad \quad \quad \quad \quad \quad \quad \quad \quad$ (c).
}}
   \caption{BGC Hidden Reflective Barriers Relation to Ornstein-Uhlenbeck Process}
   \label{Fig:BGCHiddenReflectiveBarriersRelationtoOrnsteinUhlenbeckProcess}
\flushleft
   \textbf{\footnotesize
(a). Idealized BGC Process: BGCSP involves It\^{o} processes that start at an initial position $\boldsymbol{x_0 = 0}$, where the further they drift away from the long run mean $\boldsymbol{X = 0}$, then the slower they will approach the hidden barriers $\boldsymbol{\mathfrak{B}_L}$, $\boldsymbol{\mathfrak{B}_U}$ where $\boldsymbol{x_0 > \mathfrak{B}_L}$ and $\boldsymbol{x_0 < \mathfrak{B}_U}$.\\
(b). Idealized Hybrid between BGCSP and OUP: Via the Reflection Principle, we can better see the transition between the other two paradigms.\\
(c). Idealized Ornstein-Uhlenbeck Process: OUP involves It\^{o} processes that can start at any initial position $\boldsymbol{x_0 \in \mathbb{R}}$ and will `gravitate' to the long-term mean $\boldsymbol{\alpha}$.\\
}
\end{figure}
%\FloatBarrier

\bigskip \noindent
From Figure \ref{Fig:BGCHiddenReflectiveBarriersRelationtoOrnsteinUhlenbeckProcess}, we can see that the total number of paths (as shown in blue and red) is preserved in each of the three schemes.
(\ref{Eq:BGCSP1}) is (\ref{Eq:BGC1}) repeated for convenience and (\ref{Eq:BGCSP1}) is similar to (\ref{barriers}), where $\kappa (\alpha - X_t)$ in OUP is replaced by $\kappa (\alpha -\sgn[X_t, t] \Psi (X_t, t))$ in BGCSP (with $\kappa = 1$ and $\alpha = f(X_t, t)$), which caters for a much wider set of possible paths than is possible with OPU and yet BGCSP has $\sgn[X_t, t]\Psi(X_t, t)$ instead of $X_t$.
Multiplying (\ref{Ornstein_Uhlenbeck_Process}) by $e^{\kappa t}$ and expanding gives,

\begin{eqnarray}
\begin{array}{rcl}
%------------------------------------------------------------------------------------------------------------
 e^{\kappa t} dX_{t} + \kappa e^{\kappa t} X_{t} \, dt & = & \kappa \alpha e^{\kappa t} \, dt  + \sigma  e^{\kappa t}   \, dW_{t} \nonumber \\
%------------------------------------------------------------------------------------------------------------
 d(e^{\kappa t} X_{t})  & = & \kappa \alpha e^{\kappa t} \, dt  + \sigma  e^{\kappa t}   \, dW_{t} \nonumber \\
%------------------------------------------------------------------------------------------------------------
 \displaystyle \int^{T}_{0} d(e^{\kappa t} X_{t})  & = & \displaystyle \int^{T}_{0} \kappa \alpha e^{\kappa t} \, dt  + \int^{T}_{0} \sigma  e^{\kappa t}   \, dW_{t} \nonumber \\
%------------------------------------------------------------------------------------------------------------
  e^{\kappa T} X_{T} - e^{0} X_{0}  & = & \displaystyle  \kappa \alpha \frac{e^{\kappa T} - e^0 }{\kappa}   +  \sigma \int^{T}_{0}  e^{\kappa t}   \, dW_{t} \nonumber \\
%------------------------------------------------------------------------------------------------------------
   X_{T}  & = & \displaystyle   X_0 e^{- \kappa T} + \alpha (1 - e^{-\kappa T})   +  \sigma \int^{T}_{0}  e^{-\kappa (T-t)}   \, dW_{t}.  \label{Ornstein_Uhlenbeck_Process_2}
%------------------------------------------------------------------------------------------------------------
\end{array}
\end{eqnarray}

\bigskip \noindent
Having solved the OUP SDE, we wish to determine where $X_t$ is most likely to be for $t = T$.
Taking the expectation of both sides,

\begin{equation}\label{Ornstein_Uhlenbeck_Process_3}
\begin{array}{rcl}
%------------------------------------------------------------------------------------------------------------
 \mathbb{E} \big[ X_{T} \big]  & = & \displaystyle  \mathbb{E} \Big[  X_0 e^{- \kappa T} + \alpha (1 - e^{-\kappa T})   +  \sigma \int^{T}_{0}  e^{-\kappa (T-t)}   \, dW_{t} \Big]  \\
%------------------------------------------------------------------------------------------------------------
                                    & = & \displaystyle   X_0 e^{- \kappa T} + \alpha (1 - e^{- \kappa T}),
%------------------------------------------------------------------------------------------------------------
\end{array}
\end{equation}

\bigskip \noindent
as $e^{-\kappa (T-t)}$ is deterministic.
Since $X_0 = 0$ for all BGCSPs, we will set all OUPs to start from the origin, giving,

\begin{equation}\label{Eq:E_OUP}
\mathbb{E} \big[ X_{T} \big] = \displaystyle  \alpha (1 - e^{- \kappa T}).
\end{equation}

\bigskip \noindent
We now have a stochastic argument that is a basis to justify (\ref{barriers}).
Assume that for greatest generality,

\begin{equation}\label{Eq:GreatestGenerality}
  \mathfrak{B}_L = -A(1 - e^{-\theta T}) + C, \quad \mathfrak{B}_U = A(1 - e^{-\theta T}) + C.
\end{equation}

\bigskip \noindent
We know that since $X_t$ is bi-directional, it is symmetric about the origin and so $C=0$.
By comparing (\ref{Eq:GreatestGenerality}) with (\ref{Eq:E_OUP}), we can see that,

\[
  A(1 - e^{-\theta T})  \sim \alpha (1 - e^{-\kappa T}), \quad A \sim \alpha, \quad \theta \sim \kappa ,
\]

\bigskip \noindent
where $\sim$ signifies a weak association.
Since BGCSP do not `force' the It\^{o} diffusion to the long-term mean $\alpha$, the time taken to reach $\alpha$ or $-\alpha$ under BGC would be greater than for OUP.
Hence in terms of distance,

\[
  \mathfrak{B}_L \geq -A(1-e^{-\theta T}), \quad \mathfrak{B}_U \leq A(1-e^{-\theta T}).
\]
\end{proof}

\bigskip
\begin{remark}
$A \sim \mu$ is intuitive because the greater the drift $| \mu |$, then the greater $A$, hence $|\mathfrak{B}|$, ie. $\mathfrak{B}_L$ and $\mathfrak{B}_U$.
We also notice that the diffusion term $\sigma$ doesn't contribute as much to $|\mathfrak{B}|$. 
We also know by some experimentation that $\theta \in [0,1]$.
\end{remark}

\bigskip \noindent
To extend \cite{TarantoKhan2020_1} further, the BGCSP Algorithm \ref{Alg:Algo1} is derived and simulated in the Results and Discussion section.

%-----------------------------------------------------------------------------------------
%\newpage
\bigskip
%\section{Appendix}

\begin{algorithm}\label{Alg:Algo1}
\tiny
\caption{Bi-Directional Grid Constrained (BGC) Stochastic Processes}
\# Pseudocode based on R\\
\textbf{INPUT: }\\
$\mu=drift,\text{ }\sigma=diffusion,\text{ }i=simulation\text{ }index,\text{ }$s$=\#\text{ } simulations=10,000,\text{ }t=time\text{ }steps=1001,\text{ }j=time\text{ }index,\text{ }Print\_Simulations=TRUE$\\
\textbf{OUTPUT: }\\
$ID\_value \leftarrow matrix(0:0, nrow = TimeSteps,   ncol = Simulations)$\\
$CX       \leftarrow matrix(0:0, nrow = TimeSteps,   ncol = Simulations)$\\
$T\_1000   \leftarrow matrix(0:0, nrow = Simulations, ncol = 1)$\\
$T        \leftarrow matrix(0:1000, nrow = TimeSteps, ncol = 1)$\\
$t = 1$\\
$i=1$\\
\For{$(i=1:Simulations)$} {
   $t = 1$\\
   \For{$(t=1:TimeSteps)$} {
      \uIf{(t $==$ 1)} {
         $CX[t,i]       \leftarrow 0$\\
         $ID\_value[t,i] \leftarrow 0$
      } \Else {
         $\#dt = (t/TimeSteps)$\\
         $dt = 0$\\
         $DW = rnorm(n=1, mean = 0, sd = 1)$\\
         $dW = DW$\\
         $t\_1 = t-1$\\
         $XX \leftarrow (\mu * dt + \sigma * dW)$\\
         $CX\_Now  \leftarrow CX[t\_1,i] + XX$\\
         $CX[t,i] \leftarrow CX\_Now$\\
         \uIf{$(CX\_Now > 0)$} {
            $ID\_value[t,i] \leftarrow ( CX\_Now - CX\_Now * CX\_Now / 100 )$
         } \Else {
            $ID\_value[t,i] \leftarrow ( CX\_Now + CX\_Now * CX\_Now / 100 )$
         }
      }
   }
   \uIf{(Print\_Simulations$==$TRUE)} {
      $plot(T, ID\_value[,i], type = ``l", ylim=c(yMax, yMin) )$
   } \Else {
      $lines(T, ID\_value[,i], type = ``l", ylim=c(yMax, yMin) )$
   }
   $T\_1000[i] \leftarrow sum(ID\_value[,i])$
}
\end{algorithm}
%\FloatBarrier

%%%%%%%%%%%%%%%%%%%%%%%%%%%%%%%%%%%%%%%%%%%%%
%%%%%%%%%%%%%%%%%%%%%%%%%%%%%%%%%%%%%%%%%%%%%
%%%%%%%%%%%%%%%%%%%%%%%%%%%%%%%%%%%%%%%%%%%%%
%%%%%%%%%%%%%%%%%%%%%%%%%%%%%%%%%%%%%%%%%%%%%
\bigskip
%\newpage
\section{Results and Discussion}

\begin{figure}[htb]
  \centering
  \begin{turn}{90}
  \begin{minipage}{8.5in}
  \centering
   \includegraphics[scale=0.45]{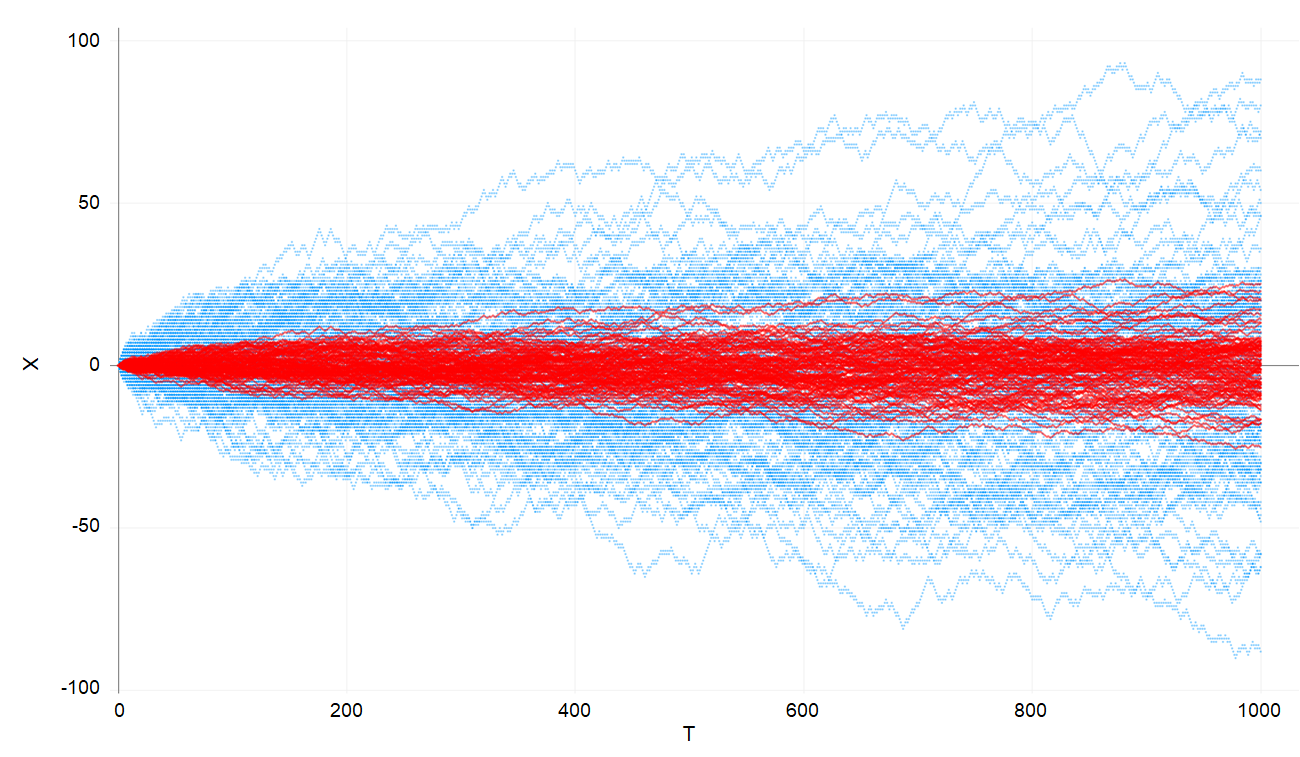}
   \includegraphics[scale=0.45]{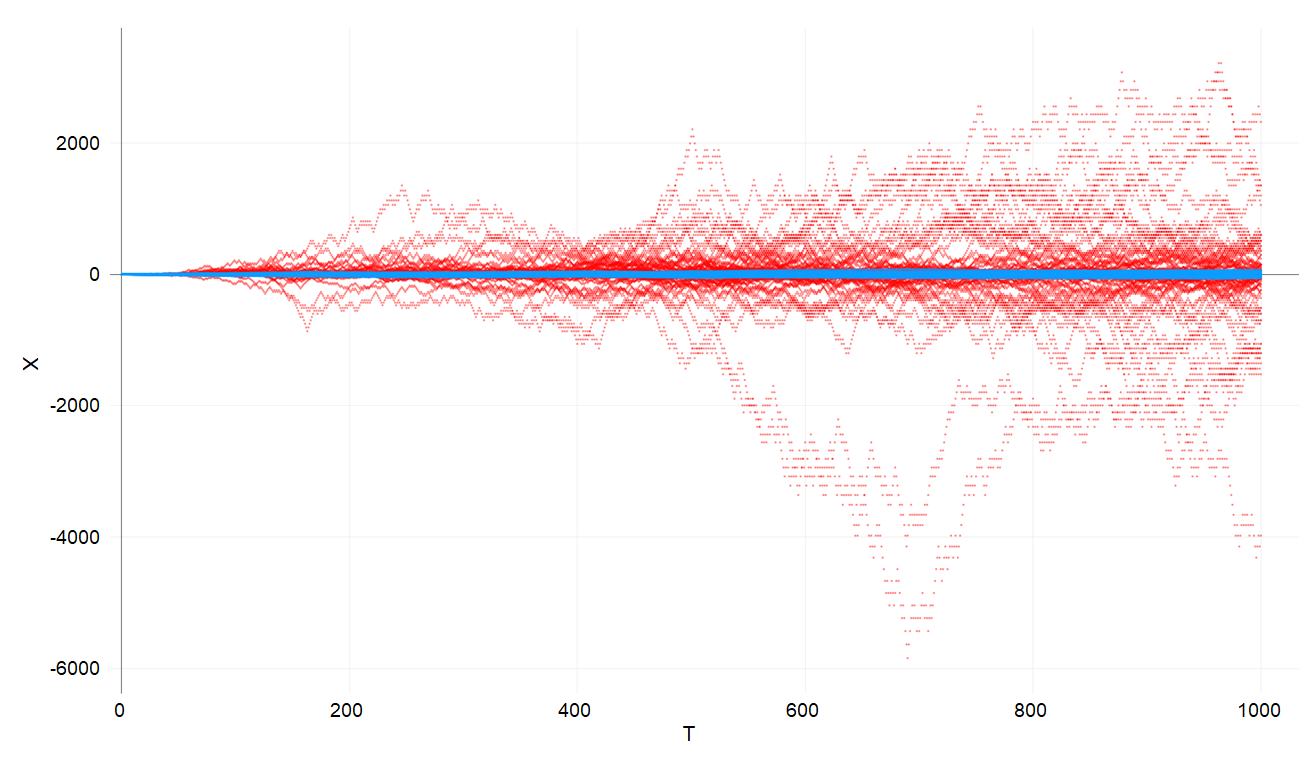}\\
   \textbf{\footnotesize \noindent{
$\quad \quad \quad \quad$ (a). $\boldsymbol{y = |x|}$ $\quad \quad \quad \quad \quad \quad \quad \quad \quad \quad \quad \quad \quad \quad \quad \quad \quad \quad \quad \quad$ (b). $\boldsymbol{(e^x + e^{-x})/\theta}$, $\boldsymbol{\theta=2000}$
}}\\
   \includegraphics[scale=0.45]{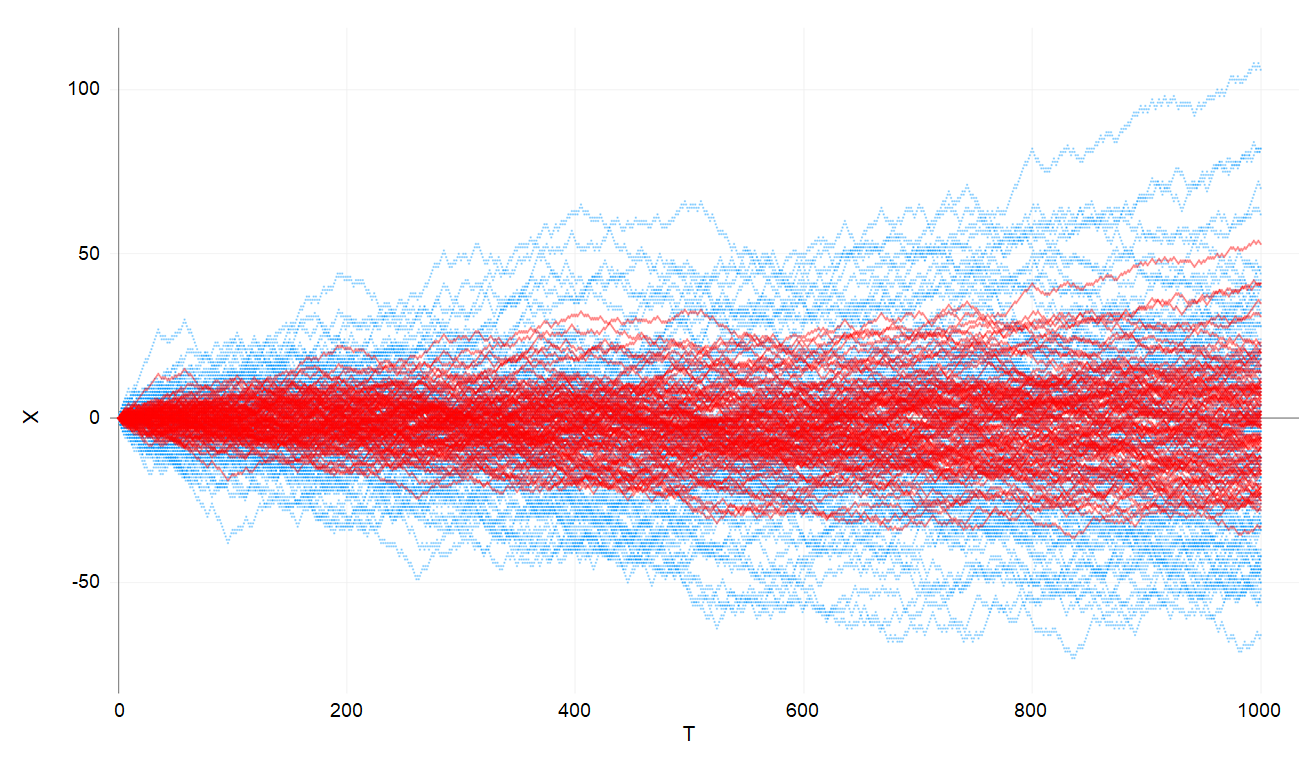}
   \includegraphics[scale=0.45]{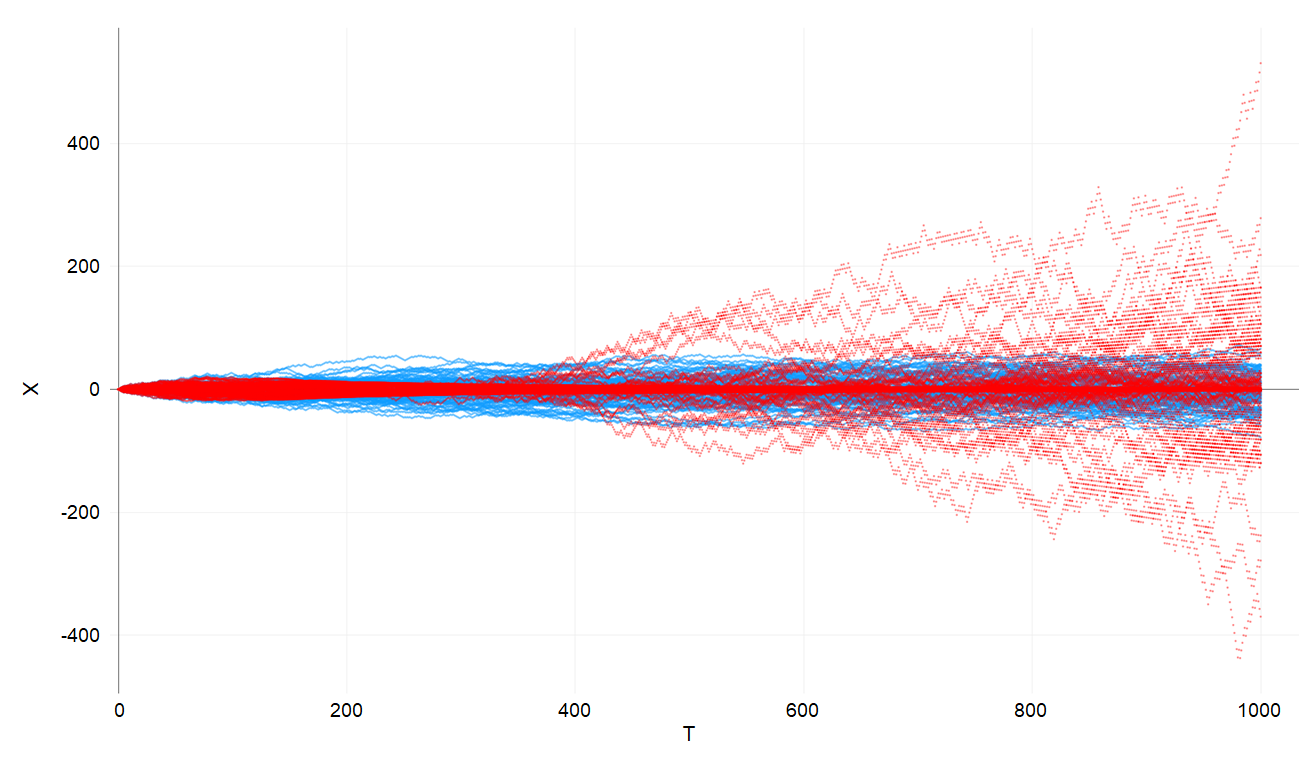}\\
   \textbf{\footnotesize \noindent{
$\quad \quad \quad \quad \quad \quad \quad $ (c). $\boldsymbol{y=x^2 t/ \theta}$, $\boldsymbol{\theta = 200}$
 $\quad \quad \quad \quad \quad \quad \quad \quad \quad \quad \quad \quad \quad \quad \quad \quad$ (d). $\boldsymbol{y=|x|^3 / \theta_1 + \theta_2}$, $\boldsymbol{\theta_1 =200}$, $\boldsymbol{\theta_2 = 5}$
}}
 \caption{Potential Convex Candidates for $\Psi (X_t, t)$}
  \label{Fig:PotentialConvexCandidatesforPsi}
%\flushleft
   \textbf{\footnotesize 
Blue = Without BGC, $\quad$ Red = With BGC\\}
\flushleft
\textbf{\footnotesize 10,000 unconstrained and 10,000 BGC It\^{o} difussions simulated.
All the potential BGC functions $\boldsymbol{\Psi(X_t, t)}$ are not suitable for the standard It\^{o} diffusions of BGCSP even though they are bi-directionally convex, except for (c).\\
(a) and (c) were rejected since they don't constrain $\boldsymbol{X_t}$ to constant hidden barriers.\\
(b) and (d) were rejected because whilst there is some initial evidence of constant hidden barriers, they become unstable over time and explode to $\boldsymbol{\pm \infty}$, faster than unconstrained It\^{o} diffusions.
}
  \end{minipage}
  \end{turn}
\end{figure}

\noindent
At this stage, we know that Use Cases I, III, IV and V -are not valid candidates for the correct type of convexity for general $\Psi (x, t)$, but for various exotic forms of It\^{o} diffusions, such as the Cox-Ingersoll-Ross (CIR) process, these use cases may be sufficient to BGC the It\^{o} process within two hidden barriers.
We also have a theoretical appreciation of what other forms of $\Psi (X_t, t)$ are valid and invalid candidates.
On the other hand, we also know that Use Case II -the parabolic cylinder -is the ideal type of convexity for general It\^{o} diffusions.
To confirm this and further eliminate any remaining use cases, we simulate them in Figure \ref{Fig:PotentialConvexCandidatesforPsi}.

\bigskip \noindent
From Figure \ref{Fig:PotentialConvexCandidatesforPsi}, in (b) and (d), we see a certain amount of constraining is occuring and that the BGC never exceeds the original simulation paths.
However, we also see in (d) that whilst there is some BGC initially, after some time the hidden barriers are unstable.

\bigskip \noindent
In (a) and (c), there is no real effective constraining since the BGCSP is now hyper-extended from the original simulation paths and so there are no hidden reflective barriers either.
When $\Psi(X_t, t) \gg X_t$ (where $X_t$ is the unconstrained Ito diffusion and $\gg$ signifies domination), then there will be a point in time where $X_t$ will flip over the origin and the contribution from the drift will be eclipsed by the $\Psi(X_t, t)$ term.
We thus have $\Psi(X_t, t)$ becoming the dominant drift term that will explode the It\^{o} diffusion beyond where the unconstrained It\^{o} diffusion would reach, away from the origin.

%\FloatBarrier

\begin{figure}[htb]
  \centering
  \begin{turn}{90}
  \begin{minipage}{8.5in}
  \centering
   \includegraphics[width=\linewidth]{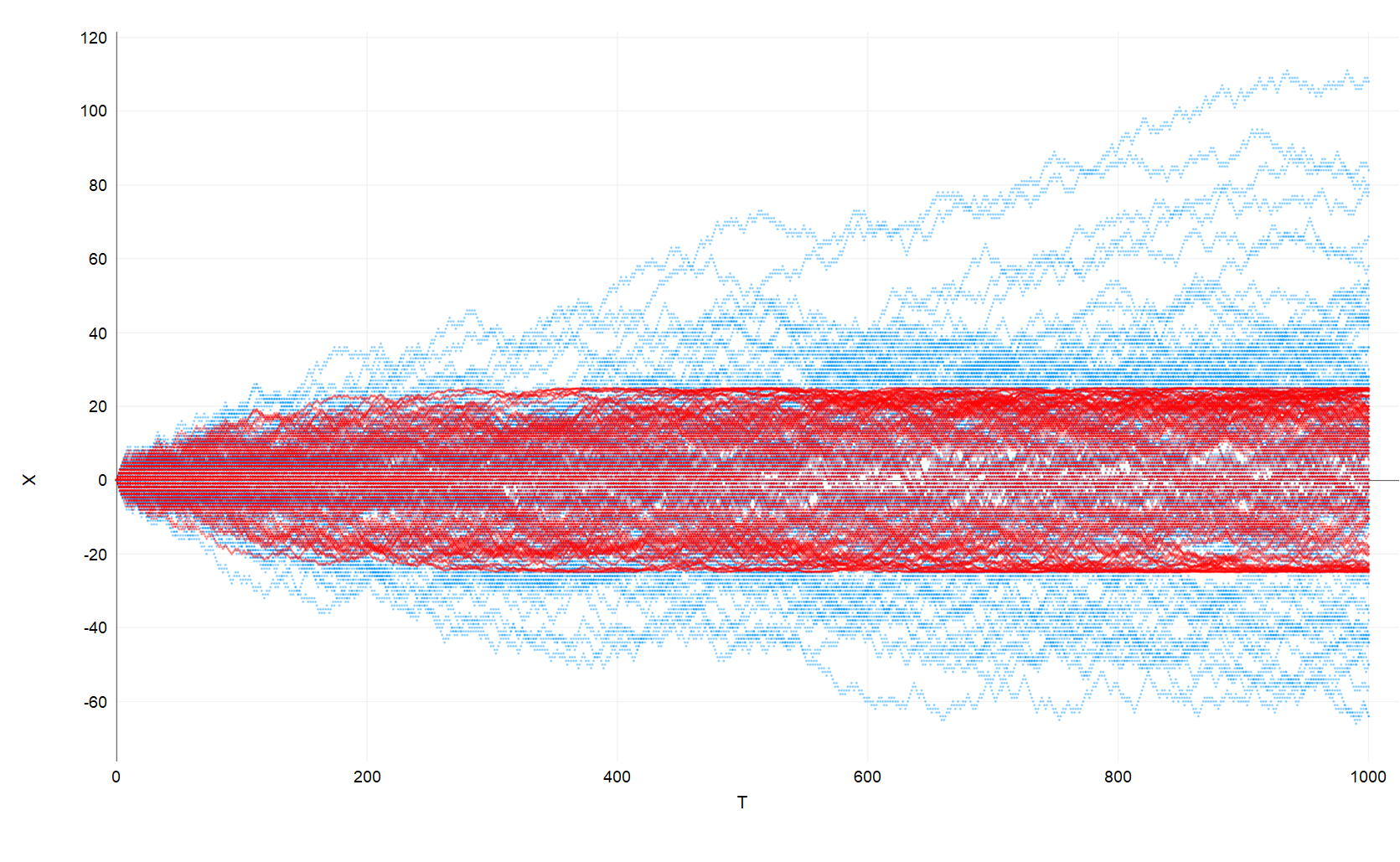}
 \caption{Ideal Selection of $\Psi (x, t)$ for BGC}
  \label{Fig:IdealSelectionofPsiforBGC}
   \textbf{\footnotesize \noindent
Blue = Unconstrained It\^{o} process, $\quad$ Red = BGC It\^{o} process.\\}
\flushleft
   \textbf{\footnotesize 10,000 unconstrained and 10,000 BGC It\^{o} difussions simulated for $\boldsymbol{\omega = 100}$ in $\boldsymbol{\Psi(X_t, t)=x^2 / \omega}$.
As more and more BGC It\^{o} processes are added to the plot, then the ideal hidden reflective barriers of BGCSP emerges, confirming that $\boldsymbol{\mathfrak{B}_L = -A (1 - e^{-\theta T})}$ and $\boldsymbol{\mathfrak{B}_U = A (1 - e^{-\theta T})}$, where $\boldsymbol{A}$ controls the position of the barrier(s) and $\boldsymbol{\theta}$ controls the speed at which the BGC It\^{o} process reaches the hidden reflective barrier.
In this plot, $\boldsymbol{A = 25}$ and $\boldsymbol{\theta = 0.01}$, but these can be altered to suit one's BGC needs.
}
  \end{minipage}
  \end{turn}
\end{figure}

\bigskip \noindent
We now examine the parabolic cylinder in far greater detail, as shown in Figure \ref{Fig:IdealSelectionofPsiforBGC}.

\begin{figure}[htb]
  \centering
  \begin{turn}{90}
  \begin{minipage}{8.5in}
  \centering
   \includegraphics[width=\linewidth]{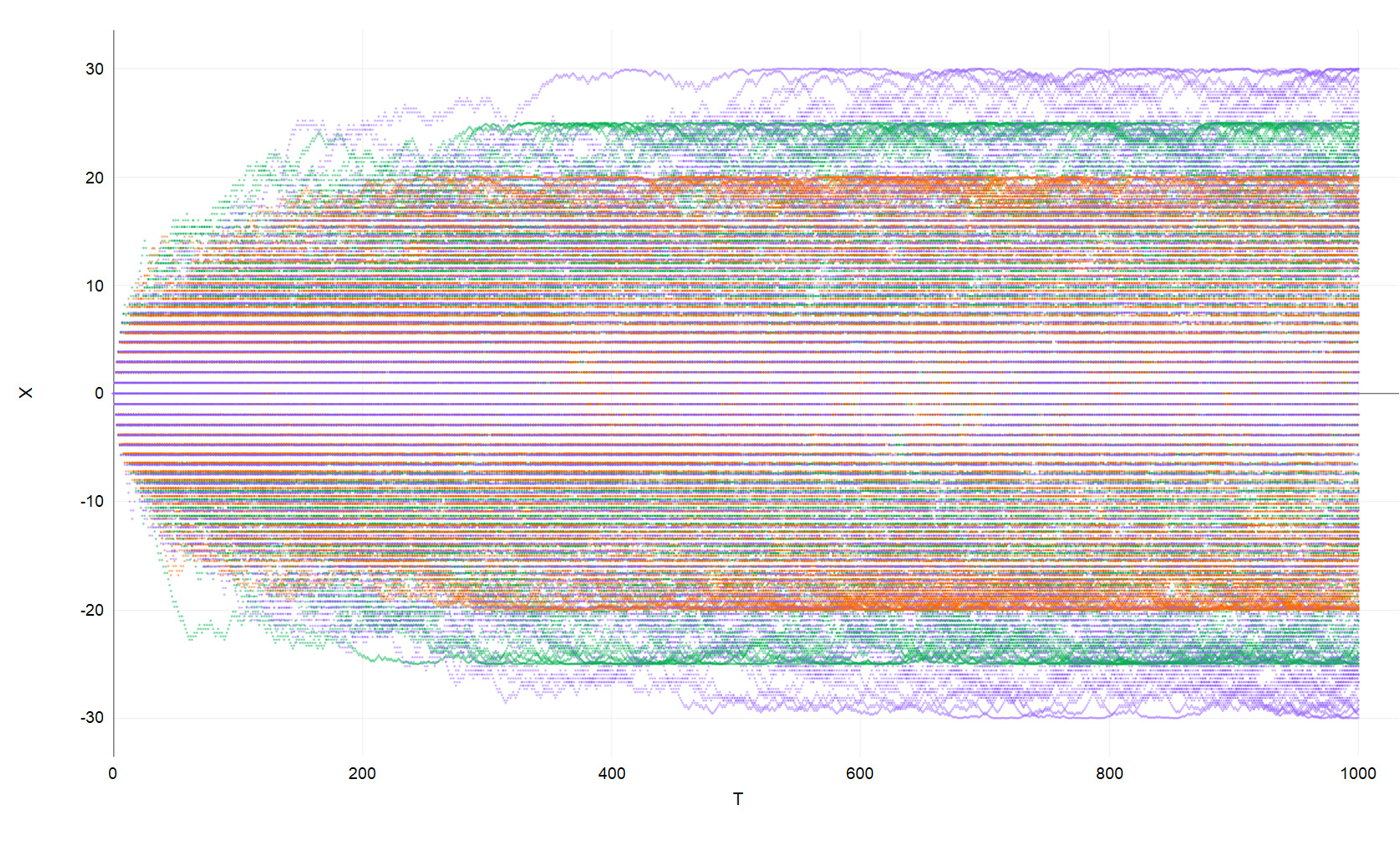}
 \caption{Parabolic Cylynder Expression of $\Psi (X_t, t)$ for Three Different Values of $A$ for $\mathfrak{B}_L$, $\mathfrak{B}_U$}
  \label{Fig:ParabolicCylynderExpressionofPsifor3DifferentValuesofA}
   \textbf{\footnotesize \noindent{
Orange = Inner-most barriers, $\quad$ Green = Middle barriers, $\quad$ Purple = Outer-most barriers.\\}}
\flushleft
   \textbf{\footnotesize 10,000 BGC It\^{o} difussions simulated for various values of $\boldsymbol{\omega}$ in $\boldsymbol{\Psi(X_t, t) = x^2 / \omega}$.
The general formula for a hidden reflective BGC barrier for a parabolic cylinder is further validated as $\boldsymbol{\mathfrak{B}_L = -A (1 - e^{-\theta T})}$ and $\boldsymbol{\mathfrak{B}_U = A (1 - e^{-\theta T})}$, where $\boldsymbol{\theta = 0.01}$ and $\boldsymbol{A=25}$ for orange, $\boldsymbol{A=30}$ for green and $\boldsymbol{A=35}$ for purple.
Note that for $\boldsymbol{\Psi (x, t) = x^2/ \omega}$, $\boldsymbol{\omega = 80}$ for orange, $\boldsymbol{\omega = 100}$ for green and $\boldsymbol{\omega = 120}$ for purple.
}
  \end{minipage}
  \end{turn}
\end{figure}

\bigskip \noindent
From Figure \ref{Fig:IdealSelectionofPsiforBGC}, we see that the hidden reflective BGC barriers can constrain the It\^{o} diffusion(s) indefinitely as it is `trapped' within the barriers.
This assumes that there are no sudden jumps (as is the case in jump-diffusion models) or changes in $X_t$ or in $\Psi (x, t)$.
To examine this in even further detail, we simulate again for different parameters, as shown in Figure \ref{Fig:ParabolicCylynderExpressionofPsifor3DifferentValuesofA}.

\bigskip \noindent
From Figure \ref{Fig:ParabolicCylynderExpressionofPsifor3DifferentValuesofA}, there is a region about the time axis where not many simulation paths visit, supporting the notion that as the paths approach the hidden barriers, they end up being `trapped' near that boundary.
Also notice how there is banding or discretization about various local times which get compressed the further they are from the origin.
The local times seem to coincide or line up most near the time axis regardless of $\omega$.

%%%%%%%%%%%%%%%%%%%%%%%%%%%%%%%%%%%%%
%\newpage
\bigskip
\section{Conclusions}

\noindent
This paper has extended the available research on BGC stochastic processes by investigating the hidden geometry of BGC functions.
The parabolic cylinder was found to be the ideal constraining mechanism for the parameter $\Psi (x, t)$, for most general unconstrained It\^{o} diffusions. 
Not any ordinary convex function will suffice and the novel `bi-directionally convex' definition was defined and adopted.
The formulas for the lower hidden reflective barrier $\mathfrak{B}_L$ and the upper hidden reflective barrier $\mathfrak{B}_U$ were derived.
This helps establish a linkage between $\Psi (X_t, t)$ of the form $x^2 / \omega$ and the resulting $\mathfrak{B}_L$ and $\mathfrak{B}_U$.
This research has applications in many fields, such as in finance where exchange rates can be constrained by `parabolic cylinder' monetary policies, such as `keep the AUD/NZD exchange rate within a range by regulating the amount of Government debt, the more it approaches the range boundaries'.
Future research in BGC can involve BGC of other important It\^{o} diffusions from other research fields and finding estimates for the distribution of the first passage time (FPT) for when $\mathfrak{B}_L$ and $\mathfrak{B}_U$ are most likely to be first hit.
We believe that there must be some mapping from $\Psi (x, t) = x^2/\omega$ to $\mathfrak{B}_L = -A (1 - e^{-\theta T})$ and to $\mathfrak{B}_U = A (1 - e^{-\theta T})$.

\FloatBarrier

\FloatBarrier

\bigskip

\bibliography{References}  % Path to your References.bib file %<--------------------------NOT OVERRIDDEN
\bibliographystyle{apalike}
%\bibliographystyle{ieeetr}
%\bibliographystyle{plain}
%\bibliographystyle{siam}
%\bibliographystyle{ieeer}
%\bibliographystyle{plainnat}
%\bibliographystyle{abbrv}

\begin{comment}
\bibliographystyle{amsplain}

\end{comment}

\end{document}